\documentstyle[11pt,amsfonts]{article}
\newtheorem{teo}{Theorem}[section]
\newtheorem{prop}[teo]{Proposition}
\newtheorem{cor}[teo]{Corollary}
\newtheorem{defi}[teo]{Definition}
\newtheorem{defi-prop}[teo]{Proposition--Definition}
\newtheorem{lem}[teo]{Lemma}
\newtheorem{rema}[teo]{Remark}

\newcommand{\Hom}{{\rm Hom}}
\newcommand{\lHom}{{\rm Hom}_{\log}}

\newcommand{\Sp}{{\rm Sp}}
\newcommand{\spe}{{\rm sp}}
\newcommand{\Spf}{{\rm Spf}}
\newcommand{\Spec}{{\rm Spec}}
\newcommand{\co}{{\cal O}}
\newcommand{\cC}{{\cal C}}

\newcommand{\cZ}{{\cal Z}}
\newcommand{\cS}{{\cal S}}
\newcommand{\cM}{{\cal M}}

\newcommand{\cJ}{{\cal J}}

\newcommand{\cA}{{\cal A}}

\newcommand{\cU}{{\cal U}}

\newcommand{\fX}{{\frak X}}
\newcommand{\fY}{{\frak Y}}
\newcommand{\fZ}{{\cal Z}}

\newcommand{\fC}{{\frak C}}

\newcommand{\fn}{{\frak n}}
\newcommand{\fm}{{\frak m}}
\newcommand{\ofX}{{\overline{\frak X}}}
\newcommand{\ofY}{{\overline{\frak Y}}}
\newcommand{\cY}{{\cal Y}}

\newcommand{\cD}{{\cal D}}

\newcommand{\cF}{{\cal F}}

\newcommand{\oh}{\overline{h}}
\newcommand{\oL}{\overline{L}}

\newcommand{\of}{\overline{f}}

\newcommand{\oOmega}{\overline{\Omega}}

\newcommand{\oK}{\overline{K}}
\newcommand{\oU}{\overline{U}}

\newcommand{\ocD}{\overline{\cal D}}

\newcommand{\ocZ}{\overline{\cal Z}}

\newcommand{\ocM}{\overline{\cal M}}

\newcommand{\ofC}{\overline{\frak C}}
\newcommand{\ocU}{\overline{\cal U}}

\newcommand{\tF}{\widetilde{F}}

\input xy
\xyoption{all}
\CompileMatrices

\begin{document}

\title{Ramification of local fields with imperfect residue 
fields }
\author{{\sc Ahmed Abbes and Takeshi Saito}}
\maketitle

\abstract{
We define two decreasing filtrations by ramification groups
on the absolute Galois group of
a complete discrete valuation field
whose residue field may not be perfect.
In the classical case where 
the residue field is perfect,
we recover the classical upper numbering filtration.
The definition uses rigid geometry 
and log-structures.
We also establish some of their properties.}

\section{Introduction}

Let $K$ be a complete discrete valuation field, and let $G$ be the Galois group
of a separable closure $\Omega$. 
Classically the ramification filtration of $G$ is defined in the case where 
the residue field of $K$ is perfect (\cite{serre1}, Chapter IV). 
In this paper, we define without any assumption on the residue field, 
two ramification decreasing filtrations 
of $G$ and study some of their properties.
Our first filtration, $(G^a)_{ a\in {\Bbb Q}_{\geq 0}}$, 
satisfies the following properties:
\begin{itemize}
\item[(i)] for a rational number $0<a\leq 1$, 
$G^a$ is the inertia subgroup of $G$;
\item[(ii)] $G^{1+}=\overline{\cup_{a>1}G^a}$ is the wild inertia 
subgroup of $G$;
\item[(iii)] it is stable by unramified extensions of $K$;
\item[(iv)]  it coincides 
with the classical upper numbering ramification filtration shifted by one, 
if the residue field of $K$ is perfect. 
\end{itemize}
When the residue field of $K$ is perfect, the classical ramification 
filtration is stable by tame base change of $K$, in the sense of (iii') 
below. In general, our ramification filtration 
(shifted by $-1$) is not stable by tame base change. 
We define a second filtration of $G$,
the logarithmic ramification filtration
$(G^a_{\log})_{a\in {\Bbb Q}_{\geq 0}}$, that has this property.
It satisfies the following:
\begin{itemize}
\item[(ii')] $G_{\log}^{0+}=\overline{\cup_{a>0}G_{\log}^a}$ 
is the wild inertia subgroup of $G$;
\item[(iii')] let $K'$ be a finite separable extension of $K$ 
of ramification index $m$, 
contained in $\Omega$. We identify the Galois group 
$G_{K'}={\rm Gal}(\Omega/K')$ 
with a subgroup of $G$. Then, we have 
$G^{ma}_{K',\log}\subset G^a_{\log}$, with equality if $K'/K$ is tamely ramified;
\item[(iv')] it coincides 
with the classical upper numbering ramification filtration, 
if the residue field of $K$ is perfect. 
\end{itemize}

We prove that both filtrations are left continuous and their jumps 
are rational. More precisely, 
for a real number $a>0$, we put $G^{a+}=\overline{\bigcup_{b>a}G^b}$
and $G^{a-}=\bigcap_{b<a}G^b$, where $b$ denotes a rational number. 
Then $G^{a-}=G^a$ if $a$ is rational, and $G^{a-}=G^{a+}$ if 
$a$ is not rational.  The same result holds true for the logarithmic 
filtration. 

In classical ramification theory, the following results play 
an important role:
\begin{itemize}
\item[1)] the graded pieces
($G^a/G^{a+}$, for $a\in {\Bbb Q}_{>1}$)
are abelian and killed by the residue 
characteristic of $\co_K$.
\item[2)] the Hasse-Arf theorem. 
\end{itemize}
Towards the first problem, we give a geometric interpretation of 
the graded pieces of our first ramification filtration. 
In a forthcoming paper, we will investigate the second problem
by comparing our theory with Kato's ramification theory for rank one 
Galois characters \cite{kato}.

In the following, we describe our constructions in more detail. 
Let $\cC$ be the Galois category of finite \'etale extensions 
of $K$, equipped with its fiber functor $\cF(L)=\Hom_K(L,\Omega)$.
The idea is to construct directly the functor $\cF^a(L)=\cF(L)/G^a$,
before defining $G^a$. 
Let $\co_K$ be the valuation ring of $K$, 
$\co$ be the integral closure of $\co_K$ in $\Omega$, and for a 
rational number $a>0$, let $\pi^a\in \co$ be an element of 
valuation ${\rm ord}_K(\pi^a)=a$. 
For a finite separable extension $L$ of $K$ with valuation ring $\co_L$, 
we have 
\[
\cF(L)=\lim_{\stackrel{\longleftarrow}{a\in {\Bbb Q}_{>0}}}
\Hom_{\co_K}(\co_L,\co/\pi^a \co).
\]
Roughly speaking, we define $\cF^a(L)$ as the set of connected components 
of $\Hom_{\co_K}(\co_L,\co/\pi^a \co)$ for some topology. 
More precisely, 
let $Z$ be a finite system of generators of $\co_L$ over $\co_K$, and 
let $I_Z$ be the kernel 
of the surjection $\co_K[X_1,\dots,X_n]\rightarrow \co_L$.
Let $D^n$ be the $n$-dimensional closed polydisc of radius one. 
For a rational number $a>0$, we define the $K$-affinoid variety 
\[
X^a_Z=D^n(|\pi|^{-a} f; \ f\in I_Z),
\]
as an affinoid sub-domain of $D^n$. The set of its $\Omega$-valued points 
is given by
\[
X^a_Z(\Omega)
=\{(x_1,\dots,x_n)\in \co^n ;\ {\rm ord}_K f(x_1,\dots,x_n)\geq a
\ {\rm for \ all}\ f\in I_Z\}.
\]
The natural map 
$X^a_Z(\Omega)\rightarrow \Hom_{\co_K}(\co_L,\co/\pi^a \co)$ is surjective,
and its fibers are discs. Hence, we define $\cF^a(L)$
as the set of geometric connected components of $X^a_Z$. 
We will show that the latter does not depend on $Z$. 
We construct the group $G^a$ from the functor $\cF^a$
by general Galois theory. 
Concretely, if $L$ is a finite Galois extension of $K$ 
contained in $\Omega$, $G(L/K)$ is the Galois group of $L/K$,  
and $(G(L/K)^a)_{a\in{\Bbb Q}_{\geq 0}}$ is the quotient filtration
of $(G^a)_{a\in{\Bbb Q}_{\geq 0}}$, then there are canonical isomorphisms 
\[
\cF^a(L)\simeq G(L/K)/G(L/K)^a.
\]
We develop also logarithmic analogues of these constructions.

\vspace{5mm}

\noindent{\it Acknowledgments.}\ 
We would like to acknowledge the deep influence of M. Raynaud 
at various stages of this work, particularly in Sections 7 and 8. 
It is a pleasure to thank him for his help and encouragements, 
for pointing out a mistake in a first 
version of this paper, and for helping to repair it. 
We would like to thank K. Fujiwara for telling us
ideas on the definition of the ramification filtration in
his unpublished work. A. Abbes wants to thank B. Poonen 
for the invitation to the University of California at Berkeley,
and T. Saito for the invitation to 
the University of Tokyo. T. Saito wants to thank 
the Max-Planck Institut, the University of Essen, and 
the University of Paris-Nord where the joint work was started, for their invitations.

\tableofcontents

\vspace{5mm}

\noindent{\it Notation.}\ 
For a finite separable extension $L$ of $K$, we denote 
by $\co_L$, ${\frak m}_L$, $\oL$, $\pi_L$ 
and $v_L$, respectively, the valuation ring, its maximal ideal, 
the residue field, a uniformizing element, and the normalized 
valuation of $L$. The integral closure of $\co_K$ in $\Omega$ is denoted 
by $\co$, and its residue field by $\oOmega$.  
We denote by $v$ the unique extension of $v_K$ to $\Omega$, 
and define an ultra-metric norm on $\Omega$ by $|x|=\theta^{v(x)}$,  
where  $0<\theta<1$ is a real number. 
We fix a uniformizer $\pi$ of $\co_K$. 
By abuse of notation, $\pi^a$ denotes an element of $\co$ 
of valuation $v(\pi^a)=a$.

For a positive integer $n$, 
the $n$-dimensional closed polydisc of radius one is denoted by $D^n$. 
For a finite set $I$ of cardinality $n$, we put $D^I=D^n$, 
where the natural coordinates are indexed by $I$. 
For a rational number $a>0$, the $n$-dimensional closed 
polydisc of radius $a$ is denoted by $D^{n,(a)}$.

\section{Filtration on Galois categories}

Let $(\cC,\cF)$ be a Galois category. Namely, there exits a profinite
group $G$ such that the fiber functor $\cF$ 
from $\cC$ to the category of finite sets is an equivalence 
between $\cC$ and the category of finite sets with continuous $G$--actions
\cite{murre}. 

\begin{prop}\label{sub-galois}
Let $\cF'$ be a functor from $\cC$ to the category of finite sets
with continuous $G$-actions, that commutes with finite disjoint unions. 
Let $\cF\rightarrow 
\cF'$ be a morphism of functors satisfying the following properties:
\begin{itemize}
\item[$1)$] $\cF(X)\rightarrow \cF'(X)$ is surjective 
for any object $X$ in $\cC$;
\item[{\rm 2)}] for any morphism $X\rightarrow Y$ such that 
$\cF(X)\rightarrow \cF(Y)$ is surjective, the following diagram 
is cocartesian
\[
\xymatrix{
\cF(X)\ar[r]\ar[d]&\cF'(X)\ar[d]\\
\cF(Y)\ar[r]& \cF'(Y)}
\]
\end{itemize}
Let $\cC'$ be the full--subcategory of $\cC$ of objects $X$ such that 
$\cF(X)\rightarrow \cF'(X)$ is bijective. 
Let $N$ be the intersection of the 
kernels of the actions of $G$ on $\cF(X)$ for $X$ in $\cC'$. 
Then, $(\cC',\cF|_{\cC'})$ is a Galois
category of group $G/N$. Moreover, for any object $X$ of $\cC$, 
we have $\cF'(X)=\cF(X)/N$. 
\end{prop}
{\it Proof.}\ The 
properties of Galois categories (cf. \cite{murre} chapter IV) for
$(\cC',\cF|_{\cC'})$ follow from those of $(\cC,\cF)$ and the following~:
\begin{itemize}
\item[a)] $\cC'$ is stable by finite  fiber-products;
\item[b)] $\cC'$ is stable by finite disjoint unions;
\item[c)] let $X\rightarrow Y$ be an effective epimorphism in $\cC$.
If $X$ is an object of $\cC'$, then $Y$ is an object of $\cC'$.
\end{itemize}
Condition b) follows from the commutativity of $\cF'$ 
with finite disjoint unions, and 
condition c) is a consequence of 2). We prove that condition 
a) is satisfied. Let $X\rightarrow Z$ and $Y\rightarrow Z$ 
be morphisms in $\cC'$. We have a commutative diagram 
\[
\xymatrix{
\cF(X\times_ZY)\ar[r]^f\ar[d]_i& \cF'(X\times_ZY)\ar[d]\\
\cF(X)\times_{\cF(Z)}\cF(Y)\ar[r]_g& \cF'(X)\times_{\cF'(Z)}\cF'(Y)}
\]
Since $i$ and $g$ are bijective and $f$ is surjective, 
then the map $f$ is bijective.  
The statement that $G/N$ is the Galois group of $\cC'$ is obvious. 

Let $X$ be an object of $\cC$, and let $Y$ be the object
given by $\cF(Y)=\cF'(X)$. Assumption 2) implies that $\cF(Y)\simeq \cF'(Y)$. 
Hence $Y$ is an object of $\cC'$, and $N$ acts trivially on $\cF(Y)=\cF'(X)$. 
We deduce a surjective map $\varphi: \cF(X)/N\rightarrow \cF'(X)$. 
Let $Z$ be the object of $\cC'$ 
given by $\cF(Z)=\cF(X)/N$, and let $\psi: \cF'(X)\rightarrow \cF'(Z)=
\cF(X)/N$ be the canonical map. We have $\psi\circ \varphi ={\rm id}$. 
So $\varphi$ is injective, and therefore bijective. 
\hfill $\Box$

\begin{rema}{\rm 
If we replace in Proposition \ref{sub-galois}
property 2) by the following weaker property:
\begin{itemize}
\item[$2')$] let $X\rightarrow Y$ be a morphism such that 
$\cF(X)\rightarrow \cF(Y)$ is surjective. If $\cF(X)\rightarrow \cF'(X)$ 
is bijective, then $\cF(Y)\rightarrow \cF'(Y)$ is bijective; 
\end{itemize}
then we still have $(\cC',\cF_{|\cC'})$ is a Galois category 
of group $G/N$, but we may have $\cF'(X)\not=\cF(X)/N$.}
\end{rema}

\begin{defi}
A filtration of the fiber functor $\cF$ is given by the following
data~:
\begin{itemize}
\item[{\rm a)}] for each rational number $v>0$, a functor $\cF^v$ 
from $\cC$ to the category of finite sets with continuous 
$G$-actions, that commutes with finite disjoint unions, and a morphism 
of functors $\varphi^v: \cF\rightarrow \cF^v$ 
satisfying properties $1)$ and $2)$ 
of Proposition $\ref{sub-galois}$;
\item[{\rm b)}] for each rational numbers $v\geq w>0$, a morphism of functors 
$\varphi_v^w:\cF^v\rightarrow \cF^w$ such that $\varphi^w=\varphi^w_v\circ 
\varphi^v$.
\end{itemize}
The filtration is said to be separated if the morphism of functors 
\begin{equation}\label{separe}
\cF\longrightarrow \lim_{\stackrel{\longleftarrow}{v\in {\Bbb Q}_{>0}}}
\cF^v
\end{equation}
is an isomorphism. 
\end{defi}

A filtration of the fiber functor $\cF$ defines a decreasing filtration by 
closed normal subgroups of the Galois group $G$, and vice versa. 
The filtration of the fiber functor is separated if and only if 
the filtration of the group is separated.

\section{The main constructions}\label{mc}

Let $\cC$ be the category of finite \'etale schemes over ${\rm Spec}(K)$.
The functor $\cF$ from $\cC$ to the category of finite sets which
maps $X$ to $\cF(X)=X(\Omega)$, makes $\cC$ into a Galois category of group $G=
{\rm Gal}(\Omega/K)$. Using the results of the previous section,
we define on $G$ two separated filtrations which generalize
the filtration by higher ramification subgroups for local fields 
with perfect residue fields.

\subsection{The non--logarithmic construction}\label{construction-1}

Let $A$ be a finite flat $\co_K$--algebra,
let $Z=(z_1,\dots,z_n)$ be a system of generators of $A$ as an 
$\co_K$-algebra,  and let $I_Z$ be the kernel of the surjection 
$\co_K[X_1,\dots,X_n]\rightarrow A$. 
For a rational number $a>0$, we define 
\[
X_Z^a=D^n(\theta^{-a}f; \ f\in I_Z).
\]
It is an affinoid sub-domain of $D^n$ since $X^a_Z=D^n(\theta^{-a}f_1,
\dots,\theta^{-a}f_m)$ for a finite system $(f_1,\dots,f_m)$ 
of generators of $I_Z$.  
The collection of these affinoid varieties 
satisfies the following properties~:
\begin{itemize}
\item[(i)] If $L=A\otimes_{\co_K}K$ is \'etale over $K$, 
by identifying a homomorphism $\varphi: L\rightarrow \Omega$ 
with the point $(\varphi(z_1),\dots,\varphi(z_n))$, 
the finite set $\cF(L)=\Hom_K(L,\Omega)$ is a subset of $X_Z^a(\Omega)$. 
\item[(ii)] If $b\geq a$ are rational numbers, then $X^b_Z$ is an affinoid 
sub-domain of $X^a_Z$. 
\item[(iii)] Let $B$ be a finite flat $\co_K$--algebra, and let 
$u: A\rightarrow B$ be a morphism of $\co_K$--algebras.
Let $Z'=(z'_1,\dots,z'_{n'})$ be a finite system of generators of $B/\co_K$,
and let $v:\{1,\dots,n\}\rightarrow \{1,\dots,n'\}$ be a map
satisfying $u(z_i)=z'_{v(i)}$. 
Then, the morphism $D^{n'}\rightarrow 
D^n$ given by $(x_1,\dots,x_{n'})\mapsto (x_{v(1)},\dots,
x_{v(n)})$ induces a morphism 
$X^a_{Z'}\rightarrow X^a_Z$. 
\end{itemize}
By (iii) above, the collection $(X^a_Z)$, where $Z$ 
runs over the finite systems of generators of $A/\co_K$,
forms a projective system. 
Let $\pi_0(X_Z^a)$ be the set of geometric connected components 
of $X_Z^a$ with respect to either the weak or the strong $G$--topology
(see \cite{bgr} 9.1.4/8). 

\begin{lem} \label{const-1}
The projective system $(\pi_0(X^a_Z))_Z$ is constant. 
\end{lem}
{\it Proof.}\ It is enough to prove that 
$\pi_0(X^a_{Z'})\simeq\pi_0(X^a_Z)$, when $Z'=(Z,t)$ and $t\in A$. 
Let $g\in \co_K[X_1,\dots,X_n]$ be a lifting of $t$,
and $A=\co_K[X_1,\dots, X_n, T]/I_{Z'}$ be the presentation 
associated with $Z'$. If $I_Z=(f_1,\dots,f_m)$, then 
$I_{Z'}=(f_1,\dots,f_m,t-g)$. Therefore, 
$X_{Z'}^a\simeq X_Z^a\times D^{1,(a)}$ and the morphism of (iii)
is the first projection. \hfill $\Box$

\vspace{0.5cm}

We put  $\displaystyle 
\cF^a(A)=\lim_{\stackrel{\longleftarrow}{Z}}\pi_0(X^a_Z)$.
 
\begin{lem}\label{extra}
{\rm i)} Let $a>0$ be a rational number, and by abuse of notation
let $\pi^a\in \co$
be an element of valuation $v_K(\pi^a)=a$. Then, the  
map $X^a_Z(\Omega)\rightarrow \Hom_{\co_K}(A,\co/\pi^a\co)$ 
sending $(x_1,\dots,x_n)$ to the homomorphism defined by 
$z_i\mapsto x_i ({\rm mod}\ \pi^a \co)$, 
is surjective and its fibers are discs. 
\item[{\rm ii)}] The natural map $\bigcup_{a>0}X^a_Z(\Omega)
\rightarrow \Hom_{\co_K}(A,\oOmega)$ induces a bijection
\[
\lim_{\stackrel{\longrightarrow} {a\in {\Bbb Q}_{>0}}}
\cF^a(A)\longrightarrow \Hom_{\co_K}(A,\oOmega).
\]
\end{lem}
{\it Proof.}\ Obvious. \hfill $\Box$

\vspace{2mm}

The definition of $X^a_Z$ may be explained more geometrically 
as follows. We put $\fX={\rm Spec}(A)$. Then the system $Z$ of generators 
defines a closed immersion $i_Z:\fX\rightarrow {\Bbb A}^n_{\co_K}$. 
In the commutative diagram 
\[
\xymatrix{
{\cF(L)=\fX(\Omega)}\ar[r]& {X^a_Z(\Omega)} \ar[d]\ar[r]& {{\Bbb A}^n(\co)}
\ar[d]\\
&{\fX(\co/\pi^a\co)}\ar[r]&{{\Bbb A}^n(\co/\pi^a\co)}}
\]
the square is Cartesian and the fibers of the vertical maps are discs.

\vspace{2mm}

Let $X={\rm Spec}(L)$ be an object in $\cC$ and $\co_L$ be the normalization
of $\co_K$ in $L$. We associate with $X$ the set $\cF^a(X)=\cF^a(\co_L)$. 
We easily check that if $X=\amalg_{i=1}^r{\rm Spec}(L_i)$, then $\cF^a(X)
=\amalg_{i=1}^r\cF^a(\co_{L_i})$. 
By (iii) above, $\cF^a$ is a functor 
from $\cC$ to the category of finite sets with continuous 
$G$-actions. By (i) and (ii) above, we 
have natural morphisms of functors $\cF\rightarrow \cF^a$
and $\cF^b\rightarrow \cF^a$, 
for rational numbers $b\geq a>0$.

\begin{teo} \label{fil1}
The projective system $(\cF^a, \cF\rightarrow \cF^a)_{a\in {\Bbb Q}_{>0}}$ 
is a separated filtration of the fiber functor $\cF$. 
\end{teo}

\begin{defi} 
The filtration $(G^a, a\in {\Bbb Q}_{\geq 0})$ 
of the Galois group $G$, induced by the filtration of 
the fiber functor $\cF$ of theorem 
$\ref{fil1}$ and extended by $G^0=G$, is called 
the ramification filtration. 
\end{defi}

\begin{rema}{\rm Let $L/K$ be a finite separable extension.
By Proposition \ref{sub-galois} and Theorem \ref{fil1}, 
$\cF^a(L)=\cF(L)/G^a$ for any rational number $a\geq 0$. 
In particular, if $L$ is a finite Galois extension of $K$ 
contained in $\Omega$, $G(L/K)$ is the Galois group of $L/K$, 
and $(G(L/K)^a,{a\in{\Bbb Q}_{\geq 0}})$ is the quotient filtration
of $(G^a,{a\in{\Bbb Q}_{\geq 0}})$, then there are canonical isomorphisms 
$\cF^a(L)\simeq G(L/K)/G(L/K)^a$, for all rational numbers $a\geq 0$.}
\end{rema}

\begin{defi}
For a real number $a>0$, we put $G^{a+}=\overline{\bigcup_{b>a} G^b}$ 
and $G^{a-}=\bigcap_{b<a} G^b$, where $b$ denotes a rational number.
\end{defi}

\begin{prop}\label{properties}
$1)$ For a rational number $0<a\leq 1$, $G^a$
is the inertia subgroup of $G$, 
and $G^{1+}$ is the wild inertia subgroup. \\[3mm] 
$2)$\ Let $K'$ be a finite separable extension of $K$ contained in $\Omega$,  
of ramification index $m$. We identify the Galois group $G_{K'}=
{\rm Gal}(\Omega/K')$ with a subgroup of $G$. Then, 
for a rational number $a>0$, we have $G^{ma}_{K'}\subset G^a$, with 
equality if $K'/K$ is unramified. \\[3mm]
$3)$\ If the residue field of $K$ is perfect, then 
the filtration $(G^a)_{a\in {\Bbb Q}_{\geq 0}}$ coincides 
with the classical upper numbering ramification filtration shifted by one.
\end{prop}

\begin{teo}\label{rational}
The ramification filtration is left continuous and its 
jumps are rational, i.e., for a real number
$a>0$, we have $G^{a-}=G^a$ if $a$ is rational,
and $G^{a-}=G^{a+}$ if $a$ is not rational.
\end{teo}

\subsection{The logarithmic construction} \label{construction-2}

Let $L$ be a finite separable extension of $K$. 
A logarithmic system of generators of $\co_L$ over $\co_K$
is a triple $(Z,I,P)$ where 
$Z=(z_i)_{i\in I}$ is a finite system of generators 
of $\co_L$ as an $\co_K$--algebra 
and $P$ is a subset of $I$ such that the set $\{z_i; i\in P\}$
contains a uniformizer 
of $\co_L$ and does not contain the zero element.
Let $(Z,I,P)$ be such a system, $e$ be the ramification 
index of $L/K$, and for $i\in P$, $e_i=v_L(z_i)$. 
Let $\co_L=\co_K[(X_i)_{i\in I}]/I_Z$ be the presentation 
associated with $Z$ and 
$(f_1,\dots,f_m)$ be a finite set of generators of $I_Z$.
For $i\in P$, we choose $g_i\in \co_K[(X_i)_{i\in I}]$ 
a lifting of the unit $u_i=z_i^{e}/\pi^{e_i}$. 
For $(i,j)\in P^2$, we choose $h_{i,j}\in \co_K[(X_i)_{i\in I}]$
a lifting of the unit $u_{i,j}=z_j^{e_i}/z_i^{e_j}$. 
Let $a>0$ be a rational number. We define 
\begin{equation}\label{log-def1}
Y^a_{Z,P}=D^I\left(
\begin{array}{clcr}
\theta^{-a}f_l\ \ \ \  \forall 1\leq l\leq m\\
\theta^{-a-e_i}(X_i^{e}-\pi^{e_i}g_i)\ \ \ \  \forall i\in P\\
\theta^{-a-e_ie_j/e}(X_j^{e_i}-X_i^{e_j}h_{i,j})\ \ \ \  
\forall (i,j)\in P^2
\end{array}
\right)
\end{equation}
as an affinoid sub-domain of $D^I$.
Obviously, this definition does not depend on the choice of the $f_l$ and 
the $g_i$. It does not depend on the choice of the $h_{i,j}$ 
because for any $(x_i)_{i\in I}\in Y^a_{Z,P}(\Omega)$ and $i\in P$, 
we have $v_L(x_i)\geq e_i$.

\begin{lem} \label{equiv-def}
{\rm 1)}\ Let $g\in \co_K[(X_i)_{i\in I}]$ and $x\in X^a_Z(\Omega)$. 
If $g\in I_Z$ then $v(g(x))\geq a$, 
and if the image of $g$ in $\co_L$ is invertible
then $g(x)$ is a unit in $\co$. \\[3mm]
{\rm 2)}\ Let $\iota\in P$ be such that $z_{\iota}$ is a uniformizer of $\co_L$,
and put $g=g_{\iota}$ and $h_i=h_{\iota,i}$ for $i\in P$. 
Then, 
\begin{equation}\label{log-def2}
Y^a_{Z,P}=D^I\left(
\begin{array}{clcr}
\theta^{-a}f_l  \ \ \ \  \forall 1\leq l\leq m\\
\theta^{-a-1}(X_\iota^e-\pi g)  \\
\theta^{-a-e_i/e}(X_i-X_\iota^{e_i}h_{i}) \ \ \ \ 
\forall i\in P
\end{array}
\right)
\end{equation}
Moreover, for any $(x_i)_{i\in I}\in Y^a_{Z,P}(\Omega)$ and $i\in P$, 
we have $v_L(x_i)=e_i$. 
\end{lem}
{\it Proof.}\ 1) is obvious. 2) We put $k_i=h_{i,\iota}$ for $i\in P$.
We may assume that $g_i=h_i^e g^{e_i}$ and $h_{i,j}=h_j^{e_i} k_i^{e_j}$
for any $(i,j)\in P^2$. We denote by $Y$ the right hand side of 
(\ref{log-def2}). Let $x=(x_i)_{i\in I}\in Y(\Omega)$. 
By 1), we have $v(g(x))=0$ and $v(h_i(x))=0$ for any $i\in P$. It follows that 
$v_L(x_i)=e_i$ for any $i\in P$.
Let $i\in P$. We have $x_i=x_\iota^{e_i}(h_i(x)+\pi^a *)$ 
(where $*$ stands for an element in $\co$).
Therefore, 
\begin{eqnarray*}
x_i^e&=&x_\iota^{e e_i}(h_i(x)+\pi^a *)^e
=\pi^{e_i}(g(x)+\pi^a*)^{e_i}(h_i(x)+\pi^a *)^e\\
&=&\pi^{e_i}(g(x)^{e_i}h_i(x)^{e}+\pi^{a}*)
=\pi^{e_i} g_i(x)+\pi^{a+e_i}*
\end{eqnarray*}
Let $(i,j)\in P^2$. Since $h_i k_i=1$ modulo $I_Z$, then 
$|h_i(x)k_i(x)-1|\leq \theta^a$. 
Thus, the relation $x_i=x_\iota^{e_i}(h_i(x)+\pi^a*)$ implies 
$x_\iota^{e_i}=x_ik_i(x)+\pi^{a+e_i/e}*=x_i(k_i(x)+\pi^a*)$. 
We deduce that 
\begin{eqnarray*}
x_j^{e_i}&=& x_\iota^{e_ie_j}(h_j(x)+\pi^a *)^{e_i}
=x_i^{e_j}(k_i(x)+\pi^a*)^{e_j}
(h_j(x)+\pi^a *)^{e_i}\\
&=& x_i^{e_j}(k_i(x)^{e_j} h_j(x)^{e_i}+\pi^a*)=
x_i^{e_j}h_{i,j}(x)+\pi^{a+e_ie_j/e}*
\end{eqnarray*}
We proved that $Y(\Omega)\subset Y^a_{Z,P}(\Omega)$.
The converse is trivial. 
\hfill $\Box$

\vspace{2mm}

The collection of  affinoid varieties $Y^a_{Z,P}$ 
satisfies the following properties~:
\begin{itemize}
\item[(i)] The finite set $\cF(L)=\Hom_K(L,\Omega)$ is canonically a subset 
of $Y_{Z,P}^a(\Omega)$. 
\item[(ii)] If $b\geq a$ are rational numbers, then $Y^b_{Z,P}$ 
is an affinoid sub-domain of $Y^a_{Z,P}$. 
\item[(iii)] Let $L'$ be a finite separable extension of $L$, 
and $(Z'=(z_i)_{i\in I'},I',P')$ 
be a logarithmic system of generators of $\co_{L'}$ over $\co_K$.  
Let $v:I\rightarrow I'$ be a map such that 
$z_{i}=z_{v(i)}$ for  $i\in I$, and $v(P)\subset P'$. 
Then, the morphism $D^{I'}
\rightarrow D^I$ given by $(x_{i})_{i\in I'}\mapsto (x_{v(i)})_{i\in I}$
induces a morphism $Y^a_{Z',P'}\rightarrow Y^a_{Z,P}$. 
\end{itemize}
Only (iii) needs a proof. Let $Z'=(z_i)_{i\in I'}$, $Z=(z_i)_{i\in I}$, 
$e_i=v_{L}(z_i)$ for $i\in P$, $e'_i=v_{L'}(z_i)$ for $i\in P'$,
$e$ be the ramification index of $L/K$, $e'$ be the ramification
index of $L'/K$, and $r=e'/e$ be the ramification index of $L'/L$. 
We have $e'_i=re_i$. 
We put $u_i=z_i^e/\pi^{e_i}$ (for $i\in P$), 
$u'_i=z_i^{re}/\pi^{re_i}$ (for $i\in P'$),
$u_{i,j}=z_j^{e_i}/z_i^{e_j}$ (for $(i,j)\in P^2$), and 
$u'_{i,j}=z_j^{re_i}/z_i^{re_j}$ (for $(i,j)\in P'^2$). 
We choose 
$\iota\in P'$ such that $z_{\iota}$ is a uniformizer of $L'$. 
We have 
\begin{eqnarray*}
u_i&=&(\frac{z_i}{z_{\iota}^{re_i}})^e (\frac{z_{\iota}^{re}}{\pi})^{e_i}
=(u'_{\iota,i})^e (u'_{\iota})^{e_i}\\
u_{i,j}&=&(\frac{z_j}{z_{\iota}^{re_j}})^{e_i}
(\frac{z_{\iota}^{re_i}}{z_i})^{e_j}=
(u'_{\iota,j})^{e_i}(u'_{i,\iota})^{e_j}.
\end{eqnarray*}
We choose $g'$, $h'_i$ and $k'_i$ (for $i\in P'$) lifting 
respectively $u'_\iota$, $u'_{\iota,i}$ and $u'_{i,\iota}$. 
Let $x'=(x_i)_{i\in I'}\in Y^a_{Z',P'}(\Omega)$, and set 
$x=(x_i)_{i\in I}$. Clearly $x\in X^a_{Z}(\Omega)$. We 
prove the other relations.  
For $i\in P$, we have $x_i=
x_{\iota}^{re_i}(h_i'(x')+\pi^a*)$. Therefore, 
\begin{eqnarray*}
x_i^{e}&=&x_{\iota}^{r e_ie }(h'_i(x')+\pi^a*)^e=\pi^{e_i}(g'(x')+\pi^a*)^{e_i}
(h_i'(x')+\pi^a*)^e\\
&=& \pi^{e_i}(g'(x')^{e_i}h'(x')^e_i +\pi^a*)=
\pi^{e_i}g'(x')^{e_i}h'_i(x')^e +\pi^{a+e_i}*
\end{eqnarray*}
Let $g_i\in \co_K[(X_j)_{j\in I}]$ be a lifting of $u_i$ (for $i\in P$). 
Since $g'^{e_i} h'^e_i$ lifts $u_i$, then 
$|(g'(x'))^{e_i}(h'_i(x'))^e
-g_i(x)|\leq \theta^a$. We deduce that 
$x_i^{e}=\pi^{e_i}g_i(x) +\pi^{a+e_i}*$. 
(In general, we cannot choose 
$g_i=g'^{e_i} h'^e_i$). The other relations are proved in 
the same way.  
\hfill $\Box$

\vspace{2mm}

By (iii) above, the collection $(Y_{Z,P}^a)$, where $(Z,I,P)$ runs over 
the logarithmic systems of generators of $\co_L$ over $\co_K$, 
forms a projective system. 
Let $\pi_0(Y^a_{Z,P})$ be the set of geometric connected components 
of $Y^a_{Z,P}$. 

\begin{lem} \label{const-2}
The projective system $(\pi_0(Y_{Z,P}^a))_{(Z,P)}$ 
is constant. 
\end{lem}
{\it Proof}.\ Let $(Z,P)$ and $(Z',P')$ be as above. By considering
$Z\amalg Z'$ and $P\amalg P'$, we may assume that $Z'=(z_i)_{i\in I'}$,
$I$ is a subset of $I'$, $Z=(z_i)_{i\in I}$ and $P\subset P'$. 
Moreover, $(Z',P')$ can be obtained from $(Z,P)$ in finitely many steps, 
each step consists in either adding an element to $I$ 
and preserving $P$, or adding the same element to $I$ and $P$. 
In the first case, the proof of Lemma \ref{const-1} shows 
that the fibers of the canonical map $Y^a_{Z',P'}\rightarrow 
Y^a_{Z,P}$ are connected. 
In the second case, we have $I'=I\amalg\{t\}$ and 
$P'=P\amalg \{t\}$. We choose $\iota \in P$ as in Lemma \ref{equiv-def}
and fix $h_t\in 
\co_K[(X_i)_{i\in I}]$ a lift of the unit $z_t/z_\iota^{e_t}$. Then
$f=X_\iota^{e_t} h_t$ lifts $z_t$. Using (\ref{log-def2}), 
we deduce that the fiber of $Y_{Z',P'}^a\rightarrow Y^a_{Z,P}$ 
above a point $x=(x_i)_{i\in I}$ is isomorphic to
\begin{eqnarray*}
\lefteqn{\{x_t\in D^1(\Omega) \ ;\ |x_t-f(x)|\leq \theta^a \ \ {\rm and}\ \ 
|x_t-x_{\iota}^{e_t}h_t(x)|\leq \theta^{a+e_t/e}\}}\\
&& \ \ \ \ \ \ \ \ \ \ \ \ =\{x_t\in D^1(\Omega) \ ;\ 
|x_t-x_{\iota}^{e_t}h_t(x)|\leq \theta^{a+e_t/e}\}, 
\end{eqnarray*}
which is connected. \hfill $\Box$

\vspace{2mm}

We put $\displaystyle 
\cF^a_{\log}(L)=
\lim_{\stackrel{\longleftarrow}{(Z,P)}}\pi_0(Y^a_{Z,P})$. 
Let $X=\amalg_{i=1}^r{\rm Spec}(L_i)$ be an object in $\cC$, where each $L_i$ 
is a finite separable extension of $K$. 
We associate with $X$ the set $\cF^a_{\log}(X)
=\amalg_{i=1}^r\cF^a_{\log}(L_i)$. 
By (iii) above, $\cF^a_{\log}$ is a functor 
from $\cC$ to the category of finite sets with continuous 
$G$-actions. By (i) and (ii) above,
we have natural morphisms of functors $\cF\rightarrow \cF^a_{\log}$
and $\cF^b_{\log}\rightarrow \cF^a_{\log}$,
for rational numbers $b\geq a>0$. 

\begin{teo} \label{fil2}
The projective system $(\cF^a_{\log}, \cF\rightarrow 
\cF^a_{\log})_{a\in {\Bbb Q}_{>0}}$ 
is a separated filtration of the fiber functor $\cF$. 
\end{teo}

\begin{defi} 
The filtration $(G^a_{\log}, a\in {\Bbb Q}_{\geq 0})$ 
of the Galois group $G={\rm Gal}(\Omega/K)$, 
induced by the filtration of the fiber functor $\cF$ of Theorem 
$\ref{fil2}$ and extended by $G^0_{\log}=G$, 
is called logarithmic ramification filtration.
\end{defi}

\begin{rema}{\rm Let $L/K$ be a finite separable extension.
By Proposition \ref{sub-galois} and Theorem \ref{fil2}, 
$\cF^a_{\log}(L)=\cF(L)/G^a_{\log}$ for any rational number $a\geq 0$. 
In particular, if $L$ is a finite Galois extension of $K$ 
contained in $\Omega$, $G(L/K)$ is the Galois group of $L/K$,  
and $(G(L/K)^a_{\log},{a\in{\Bbb Q}_{\geq 0}})$ is the quotient filtration
of $(G^a_{\log},{a\in{\Bbb Q}_{\geq 0}})$, then there
are canonical isomorphisms 
$\cF^a_{\log}(L)\simeq G(L/K)/G(L/K)^a_{\log}$, 
for all rational numbers $a\geq 0$.}
\end{rema}

\begin{defi}
For a real number $a\geq 0$, we put $G^{a+}_{\log}
=\overline{\bigcup_{b>a} G^b_{\log}}$,  
and if $a>0$, we put $G^{a-}_{\log}=\bigcap_{b<a} G^b_{\log}$, 
where $b$ denotes a rational number. 
\end{defi}

\begin{prop}\label{log-properties}
$1)$ For a rational number $a>0$, we have 
$G^{a+1}\subset G^a_{\log}\subset G^a$.\\[3mm]
$2)$ $G^{0+}$ is the wild inertia subgroup of $G$.\\[3mm]
$3)$\ Let $K'$ be a finite separable extension of $K$ contained in $\Omega$,  
of ramification index $m$. We identify the Galois group $G_{K'}=
{\rm Gal}(\Omega/K')$ with a subgroup of $G$. Then, 
for a rational number $a>0$, we have 
$G^{ma}_{K',\log}\subset G^a_{\log}$, with 
equality if $K'/K$ is tamely ramified. \\[3mm]
$4)$\ If the residue field of $K$ is perfect, then 
the filtration $(G^a_{\log})_{a\in {\Bbb Q}_{\geq 0}}$ coincides 
with the classical upper numbering ramification filtration.
\end{prop}

\begin{teo}\label{log-rational}
The logarithmic ramification filtration is left continuous and its 
jumps are rational, i.e., for a real number
$a>0$, we have $G^{a-}_{\log}=G^a_{\log}$ if $a$ is rational,
and $G^{a-}_{\log}=G^{a+}_{\log}$ if $a$ is not rational.
\end{teo}

As in the non-logarithmic case, we have a geometric interpretation 
of $Y^a_{Z,P}$. For this purpose, 
we fix the following logarithmic structures. 
Let $L$ be a finite separable extension of $K$. 
The valuation ring $\co_L$ has a canonical log-structure 
$(M_L,\alpha_L)$ given by the multiplicative monoid $M_L=\co_L-\{0\}$
and the natural morphism of monoids $\alpha_L:\co_L-\{0\}\rightarrow \co_L$.
There is a canonical morphism of log-structures  $(\co_K,M_K)
\rightarrow (\co_L,M_L)$. 
Let $a>0$ be a rational number. 
The group $\{y\in \co \ ;\ v(y-1)\geq a\}$ 
acts on the monoid $\co-\{0\}$. The quotient
\[
M^a=(\co-\{0\})/\{y\in \co \ ;\ v(y-1)\geq a\}.
\]
has a canonical monoid structure. 
The morphism of monoids $\co-\{0\}\rightarrow 
\co$ induces a morphism of monoids 
$\alpha^a:M^a\rightarrow \co/\pi^a \co$. 
We have $M^a=\sqcup_{v\geq 0}M^a(v)$ 
where $M^a(v)=\{x\in \co \ ;\ v(x)=v\}/
\{y\in \co\ ;\ v(y-1)\geq a\}$. 
For a rational number $v\geq 0$, the  map
\[
M^a(v)\longrightarrow 
\{x\in \co\ ;\ v(x)=v\}/\pi^{a+v}\co,
\]
defined by $x\mapsto x$, is a well defined bijection. From this 
description, we see that 
$(M^a,\alpha^a)$ defines a log-structure on $\co/\pi^a\co$, i.e. 
$\alpha^a$ induces an isomorphism 
$M^a(0)=(\alpha^a)^{-1}((\co/\pi^a\co)^*)\rightarrow (\co/\pi^a\co)^*$. 
There is a natural morphism of log-structures $(\co_K,M_K)\rightarrow
(\co/\pi^a\co, M^a)$. 

Let $\lHom(\co_L,\co/\pi^a\co)$ be the set 
of morphisms $(\co_L,M_L)\rightarrow (\co/\pi^a\co,
M^a)$ of $(\co_K,M_K)$-log-structures. Such a morphism 
is a pair $(f,g)$ where 
$f:\co_L\rightarrow \co/\pi^a\co$ is a morphism of $\co_K$--algebras
and $g:M_L\rightarrow M^a$ is a morphism of monoids over $M_K$, such that 
$f\circ \alpha_L=\alpha^a\circ g$. 

\begin{lem} \label{hom-log}
Let $\pi_L$ be a uniformizer of $L$, $e$ be the 
ramification index of $L/K$,  
$u\in \co_L$ be the unit such that $\pi_L^e=u\pi$. 
There exists a well defined map
\[
\begin{array}{clcr}
\lHom(\co_L,\co/\pi^a\co)&\longrightarrow& \Hom(\co_L,\co/\pi^a\co)
\times M^a(\frac 1 e)\\
(f,g)&\longrightarrow& (f,g(\pi_L))
\end{array}
\]
which induces a bijection between $\lHom(\co_L,\co/\pi^a\co)$ and the set 
of pairs $(f,t)$ satisfying the equations 
\[
v(t^e-\pi f(u))\geq a+1\ \ \ \ {\rm and}\ \ \ \ 
v(t-f(\pi_L))\geq a.
\]
\end{lem}
{\it Proof.}\ Observe that $g$ is determined by its restriction 
to $\co_L^*$ and by $g(\pi_L)$. For $x\in \co_L^*$, 
$f(x)\in (\co/\pi^a\co)^*$
and $g(x)=f(x)$ via the isomorphism  
$\alpha^a:M^a(0)\simeq (\co/\pi^a\co)^*$. 
It follows that $(f,g)$ is completely determined by 
$f$ and $g(\pi_L)\in M^a$. Since $g(\pi)=\overline{\pi}\in M^a(1)$, 
then $g(\pi_L)^e\in M^a(1)$, so $g(\pi_L)\in M^a(\frac 1 e)$. 
Hence, the map in the Lemma is a well defined injection. 
Conversely, a pair $(f,t)$ as in the Lemma defines a morphism of monoids 
$g:M_L\rightarrow M^a$ over $M_K$, by $g(x\pi_L^n)=
(\alpha^a)^{-1}(f(x)).t^n$ where $x\in \co_L^*$. 
It is clear that $(f,g)\in \lHom(\co_L,\co/\pi^a\co)$. 
\hfill $\Box$

\vspace{2mm}

Let $(Z,I,P)$ be a logarithmic 
system of generators of $\co_L/\co_K$.
By Lemma \ref{hom-log}, there exists a surjective map 
\[
Y_{Z,P}^a(\Omega)\longrightarrow \lHom(\co_L,\co/\pi^a\co)
\]
with connected fibers. 
More geometrically, we may interpret our construction 
as follows. Let $\fX={\rm Spec}(\co_L)$ and let 
\[
{\cal A}={\rm Spec}(\co_K
[X_i (i\in I), 
(\frac{X_i^{e}}{\pi^{e_i}})^{\pm1} (i\in P),
(\frac{X_i^{e_j}}{X_j^{e_i}})^{\pm1}
((i,j)\in P^2)])
\]
with the log structure defined by 
${\Bbb N}\to \co_{\cal A}$, $1\mapsto X_{\iota}$, 
where $z_\iota$ is a uniformizer of $\co_L$.
We have an exact closed immersion $i:\fX\to {\cal A}$.
Then, the diagram
\[
\xymatrix{
Y^a_{Z,P}(\Omega)\ar[r]\ar[d]& {{\cal A}(\co)^{\log}}\ar[d]\\
\fX(\co/\pi^a\co)^{\log}\ar[r]& {{\cal A}(\co/\pi^a\co)^{\log}}}
\]
is Cartesian, where the superscript $\log$ indicates
the valued points as log schemes.
The horizontal maps are injective
and the fibers of the vertical maps are products of discs.

\section{Normalized integral models of affinoid varieties}\label{nim}

In Raynaud's theory \cite{bl1}, (quasi-compact and quasi-separated)
rigid $K$-spaces are the generic fibers 
of (quasi-compact) admissible $\co_K$-formal schemes.
Concretely, let $A_K$ be an affinoid $K$-algebra, and  
let $A$ be an $\co_K$--algebra which is $\pi$--adically complete and 
topologically of finite type over $\co_K$. If $A\otimes_{\co_K}K=A_K$, then 
the formal scheme ${\rm Spf}(A)$ is a model of the affinoid 
variety ${\rm Sp}(A_K)$.  
To construct such a model, we can start from a 
surjective map $\rho: K\langle X_1,\dots,X_n\rangle \rightarrow A_K$, 
and take $A_\rho=\rho(\co_K\langle X_1,\dots,X_n\rangle )$.

\begin{lem}\label{red-sup}
Assume that $A_K$ is reduced. Then,
\begin{itemize}
\item[{\rm i)}] $A_\rho\subset A=\{f\in A_K; |f|_{\sup}\leq 1\}$ 
and $A$ is the integral closure of $A_\rho$ in $A_K$. 
\item[{\rm ii)}] If $A_\rho\otimes_{\co_K} \oK$ is reduced then $A=A_\rho$.
\end{itemize}
\end{lem}
{\it Proof.}\ i) follows from \cite{bgr} 6.3.4/1. 
ii) Let $|\ |_\rho$ be the residue norm on $A_K$ defined by $\rho$. 
Then $A_\rho$ is the unit ball in $A_K$
for $|\ |_\rho$. We have $|\ |_{\sup}\leq 
|\ |_\rho$, and 
$|\ |_{\sup}=|\ |_\rho$ if and only if 
$A_\rho\otimes_{\co_K} \oK$ is reduced 
(\cite{blr4} proposition 1.1). \hfill $\Box$

\vspace{2mm}

Since $\co_K$ is a discrete valuation ring, then the unit ball $A\subset A_K$
is topologically of finite type over $\co_K$
although $A\otimes_{\co_K} \oK$ might not be reduced. However, 
we have the following finiteness theorem of Grauert and Remmert
(\cite{blr4} theorem 1.3; see also \cite{gr} and \cite{bgr} 6.4.1/3):

\begin{teo}\label{ftgr}
Let $A_K$ be a geometrically reduced affinoid $K$-algebra. Then,
there exists a finite 
separable extension $K'$ of $K$ such that the unit ball 
$A'\subset A_K\otimes_K K'$ has a geometrically reduced special 
fiber $A'\otimes_{\co_{K'}}\oK'$. 
Moreover, the formation of $A'$ commutes with any finite extension of $K'$.
\end{teo}

Let $A_K$ be a geometrically reduced affinoid $K$-algebra. 
We think of the collection of $\co_{K'}$-formal 
schemes ${\rm Spf}(A')$, where $K'$ and $A'$ are as in Theorem $\ref{ftgr}$,
as a unique model of ${\rm Sp}(A_K)$ over $\co$.
We call it the normalized integral model of ${\rm Sp}(A_K)$ over $\co$.
We say that the normalized integral model of ${\rm Sp}(A_K)$ is defined 
over $\co_{K'}$ if the unit ball $A'\subset A_{K'}$ has a geometrically 
reduced special fiber $A'\otimes_{\co_{K'}}\oK'$.

\begin{prop}\label{l-3}
Let $\fX$ be a quasi-compact
$\co_K$-formal scheme that is flat and locally of topological
finite type over $\co_K$. Assume that the closed fiber $\fX_s$ 
and the generic fiber $\fX_\eta$ (as a rigid-space) are geometrically 
reduced. Then the sets of their geometric connected 
components are isomorphic. 
\end{prop}
{\it Proof.}\ It is enough to prove that the sets of connected 
components of $\fX_s$ and $\fX_\eta$ are isomorphic. 
The set of connected components of $\fX_s$ is equal to the set 
of connected components of $\fX$. We  assume that $\fX$ is connected 
and we prove that $\fX_\eta$ is connected. 
Let $\fX=\cup_{i=1}^n {\rm Spf}(A_i)$ be a finite open covering of $\fX$
by connected affine formal schemes, 
and let $\fX_\eta=\cup_{i=1}^n {\rm Sp}(A_i\otimes_{\co_K} K)$ 
be the induced admissible covering of $\fX_\eta$. 
By Lemma \ref{red-sup} ii), $A_i$ is the unit ball for 
the sup norm on  $A_i\otimes_{\co_K} K$. If $e\in
A_i\otimes_{\co_K} K$ is a non trivial idempotent, 
then $|e|_{\sup}=1$ and $e\in A_i$.
We deduce that ${\rm Sp}(A_i\otimes_{\co_K} K)$ is connected.  

Assume that $\fX_\eta$ is not connected. Then there exists 
a partition $\{1,\dots,n\}=I\amalg J$ such that 
\[
\bigcup_{i\in I}{\rm Sp}(A_i\otimes_{\co_K}K)\cap 
\bigcup_{j\in J}{\rm Sp}(A_j\otimes_{\co_K}K)=\emptyset.
\]
We consider the open formal subschemes 
$\fX'=\bigcup_{i\in I}{\rm Spf}(A_i)$ and 
$\fX''=\bigcup_{j\in J}{\rm Spf}(A_j)$ 
of $\fX$. Since $\fX$ is connected
then there is $x\in \fX'_s\cap \fX''_s$ which is closed in $\fX_s$.
Let $\spe:\fX_\eta\rightarrow \fX_s$ be the specialization map.  
Since $\fX'$ and $\fX''$ are open, then $\spe^{-1}(x)\subset \fX'_\eta$
and $\spe^{-1}(x)\subset \fX''_\eta$. We get a contradiction 
because $\spe^{-1}(x)\not=\emptyset$ by \cite{bl1} proposition 3.5.
(In fact $\spe^{-1}(x)$ is geometrically connected
by Lemma \ref{red-sup} ii) and \cite{bosch} satz 6.1.) 
\hfill $\Box$

\begin{cor}
Let $X$ be a geometrically reduced 
affinoid variety over $K$, $\fX$ be its normalized integral 
model over $\co$, and $\ofX$ be its special fiber.
Then the sets of geometric connected components of $X$ and 
$\ofX$ are isomorphic. 
\end{cor}

\begin{prop}\label{3.5}
Let $X$ and $Y$ be geometrically reduced affinoid varieties over $K$,
and let $f:X\rightarrow Y$ be a finite flat morphism of degree $p$. 
Let $\fX$ and $\fY$ be their normalized integral models over $\co$,
and let $f:\fX\rightarrow \fY$ be the canonical extension of $f$ 
(which is finite by the finiteness theorem of Grauert and Remmert).
Let $K'$ be a finite separable extension of $K$ such that $\fX$ 
and $\fY$ are defined over $\co_{K'}$. Assume that 
\begin{itemize}
\item[{\rm i)}] $f:\fX\rightarrow \fY$ is flat;
\item[{\rm ii)}] the special fibers $\ofX$ and $\ofY$ are 
(geometrically) reduced and irreducible, and $\ofY$ is normal; 
\item[{\rm iii)}] there exits a function $\tau$ on $\fX$
that generates $\co_{X}$ over $\co_{Y}$, i.e. $\co_X=\co_Y[\tau]$.
\end{itemize}
Then, there exist a function $\xi$ on $\fY$ and an element $\beta\in\co_{K'}$
such that $(\tau+\xi)/\beta\in \co_{\fX}$, 
and $\fX\simeq \Spf(\co_{\fY}[(\tau+\xi)/\beta])$.  
\end{prop}
{\it Proof.}\ Let $\cM$ be a locally free $\co_{\fY}$-module of finite type,
and let $m\in \Gamma(\fY,\cM)$ be a non-zero section. Let 
$\cJ(m)\subset \co_{\fY}$ be the ideal of coordinates of $m$,
i.e. the ideal of $\fY$ generated locally by the coordinates of $m$ in a basis 
of $\cM$ over $\co_{\fY}$. The closed fiber of $\fY\rightarrow \Spf(\co_{K'})$
is geometrically reduced and irreducible, and $\cJ(m)$ does not vanish 
identically on the rigid fiber. Then we can apply the maximum principle
of \cite{blr3} proposition 5.2 to $\cJ(m)$ and 
the morphism $\fY\rightarrow \Spf(\co_{K'})$. So there exists 
$\beta\in \co_{K'}$ such that $\cJ(m)\subset \beta\co_{\fY}$ and 
the open part of $\fY$ where $\cJ(m)=\beta\co_{\fY}$ is not empty. 
We deduce that $m/\beta\in \Gamma(\fY,\cM)$ and its residue class 
in $\ocM=\cM\otimes_{\co_{\fY}}\co_{\ofY}$ does not vanish. 

We apply this construction to $\cM=\co_{\fX}/\co_{\fY}$
and $m$ the residue class of $\tau$.  
Then, there exist $\beta\in \co_{K'}$ and $\xi\in \co_{\fY}$ such that 
$\gamma=(\tau+\xi)/\beta\in \co_{\fX}$ and $\overline{\gamma}\not\in 
\co_{\ofY}$. Since $\co_\fX$ is locally free of rank $p$ over $\fY$, 
then $\gamma$ annihilates its characteristic polynomial $F(T)\in \co_\fY[T]$,
which is an integral polynomial of degree $p$. We claim that the natural
surjection $\co_\fY[T]/F\rightarrow \co_\fY[\gamma]$ is bijective.
Since $\co_\fY[T]/F$ and $\co_\fY[\gamma]$ are flat over $\co_{K'}$,
it is enough to show that this map is bijective after tensorization
with $K'$. The latter is true because both $\co_Y[T]/F$ 
and $\co_X=\co_Y[\gamma]$ are finite and flat of rank $p$ over $Y$. 
We consider the natural maps
\[
\fX\longrightarrow \cZ=\Spf(\co_{\fY}[\gamma])
\longrightarrow \fY, 
\]
and we denote by $x$ and $y$ 
the generic points of the special fibers of respectively 
$\fX$ and $\fY$.  
Since $\overline{\gamma}\not\in \co_{\ofY}$, 
$\overline{\gamma}$ is integral 
over $\co_{\ofY}$, and $\ofY$ is normal, 
then $\overline{\gamma}\not\in \kappa(y)$.
Thus $\kappa(x)=\kappa(y)[\overline{\gamma}]$ because
$[\kappa(x):\kappa(y)]=p$.
Hence $\fZ$ has a unique point $z$ above $y$ and 
$\kappa(z)=\kappa(x)$. The claim above implies that $\fZ$ 
is finite and flat of rank $p$ over $\fY$. We deduce that
the closed fiber of $\cZ$ is reduced and irreducible. 
Thus $\fX=\fZ$ by Lemma \ref{red-sup}. \hfill $\Box$

\vspace{2mm}

Let $X$ be a geometrically reduced 
affinoid variety over $K$, $\fX$ be its normalized integral 
model over $\co$, and $\ofX$ be its special fiber. 
There is a natural $\co$-semi-linear action of $G={\rm Gal}(\Omega/K)$ 
on $\fX$, which induces an $\oOmega$-semi-linear action on $\ofX$,
called the {\em geometric monodromy}.  
More precisely, let $K'$ be a finite Galois extension of 
$K$ over which the normalized integral model of $X$ is defined; we denote 
it by $\fX_{\co_{K'}}$. The natural $K'$-semi-linear action 
of $G$ on $X_{K'}$ extends to an $\co_{K'}$-semi-linear action 
on $\fX_{\co_{K'}}$. If $K''$ is a finite Galois extension of $K$ 
containing $K'$, then $\fX_{\co_{K''}}=\fX_{\co_{K'}}\times_{\co_{K'}}
\co_{K''}$ and the semi-linear actions of $G$ on the left and right 
hand sides are compatible. 
The restriction to the inertia  
of the action of $G$ on $\ofX$ is $\oOmega$-linear.

\section{Continuity of connected components}

Let $A$ be a geometrically reduced 
affinoid $K$-algebra, $X=\Sp(A)$ be the associated affinoid 
variety,  and $f_1,\dots,f_m\in A$ be holomorphic functions.  
For every rational number $r$, 
we consider the rational sub-domain $X^r=X(\theta^{-r}f_1,\dots,
\theta^{-r}f_m)$, and denote
by $\pi_0(X^r)$ the set of its geometric connected components.
The goal of this section is to prove the following:

\begin{teo}\label{cont}
For any rational number $a$,  
there exist rational numbers $a_1,\dots, a_{n+1}$
with $-\infty=a_0 < a_1<\dots <a_n<a_{n+1}=a$, 
such that the cardinality 
of $\pi_0(X^r)$ remains constant for all rational numbers $r\in ]a_i,a_{i+1}]$.
\end{teo}
By the maximum modulus principle
we may assume $|f_i|_{\sup}\leq 1$ for all $i$. 
Let $e$ be an integer, and $C=\Sp(K\langle z,t\rangle/zt-\pi^e)$
be the annulus $\{z \in \Omega ; 0 \leq v(z)\leq e\}$. 
We consider the following relative situation
\[
\varphi:Y=\Sp\left(A\langle z,t,u_1,\dots,u_m\rangle/(f_i-u_iz, zt-\pi^e)
\right)
\longrightarrow C.
\]  
Then $\varphi$ parameterize the family of affinoid varieties $X^r$ for 
$0\leq r\leq e$. More precisely, we have $Y_z\simeq X^{v(z)}$
for $z\in C(\Omega)$. 
For rational numbers $0<s\leq r $ and $a\in \Omega$, we put
\begin{eqnarray*}
C^+(s,r)&=&\{x\in \Omega; \ s\leq v(x)\leq r \},\\
C^-(s,r)&=&\{x\in \Omega; \ s < v(x) < r \},\\
B^-(a,r)&=&\{x\in \Omega; \ v(x-a)> r \}.
\end{eqnarray*}

\begin{prop}\label{p-3}
After a finite separable base change $K'$ of $K$, 
there exists an admissible formal model $\fC'$
of $C_{K'}$ over $\co_{K'}$  such that 
if $\ofC'$ denotes the special fiber of $\fC'$
and $\spe:C_{K'}=\fC'_{K'}\rightarrow \ofC'$ denotes the specialization map, 
then the following properties hold:
\begin{itemize}
\item[{\rm i)}] $\ofC'$ is reduced, and for any closed point $x$ 
of $\ofC'$, either $\spe^{-1}(x)=C^-(s,r)$, or 
$\spe^{-1}(x)\subset B^-(a,r)\subset C_{K'}(\Omega)$;
\item[{\rm ii)}] there exists a covering $\ofC'=\cup_{i=1}^n U_i$ 
by locally closed subschemes of $\ofC'$ such that the cardinality of 
$\pi_0(Y_x)$ is constant for $x\in \spe^{-1}(U_i)$. 
\end{itemize}
\end{prop}
{\it Proof.}\ Let $\cA\subset A$ be the unit ball for the sup-norm
(which is $\pi$-adically complete and topologically of finite 
type over $\co_K$). We consider the following model of the map $\varphi$:
\[
\phi:\cY=\Spf\left(\cA\langle z,t,u_1,\dots,u_m\rangle/(f_i-u_iz, 
zt-\pi^e)\right)
\longrightarrow \fC=\Spf(\co_K\langle z,t\rangle/zt-\pi^e).
\]  
By the main theorem of \cite{blr4}, there exists a diagram 
\[
\xymatrix{
{\cZ''}\ar[d]&&\\
\cY''\ar[r]\ar[d]&\cY'\ar[d]\ar[r]& \cY\ar[d]\\
\fC''\ar[r]^h&\fC'\ar[r]&\fC}
\]
where the squares are Cartesian and 
\begin{itemize}
\item[1)] $\fC'\rightarrow \fC$ is an admissible blow-up;
\item[2)] $\fC''\rightarrow \fC'$ is quasi-finite, flat and surjective,
and moreover it is \'etale over the generic fibers;
\item[3)] $\cZ''\rightarrow \cY''$ is finite and induces an isomorphism 
over the generic fibers;
\item[4)] $\cZ''\rightarrow \fC''$ is flat and has geometrically reduced
fibers.
\end{itemize}
Properties 2), 3) and 4) are preserved after 
any base change $\cS'\rightarrow \fC'$. By semi-stable reduction 
theorem for curves, there exists a finite 
separable extension $K'$ of $K$ 
and a morphism $\cS'\rightarrow \fC'_{\co_{K'}}$,
such that the induced morphism on the generic fibers is an isomorphism
and  $\cS'$ satisfies property {\rm i)} of Proposition \ref{p-3}.
In the sequel we replace $\fC'$ by $\cS'$.

There exists a stratification $\ofC''=\amalg_i V_i$ by locally 
closed connected subschemes of $\ofC''$ such that 
$\#\pi_0(\ocZ''_x)$ is constant for $x\in V_i$.
Thus by Proposition \ref{l-3}, $\#\pi_0(Y''_x)$ is constant for $x\in 
\spe^{-1}(V_i)$. Let $\oh:\ofC''\rightarrow \ofC'$ 
be the map induced by $h$ on the special fibers. 
We put $U_i=\oh(V_i)$. If $V_i$ is an open subscheme of $\ofC''$, 
then $U_i$ is open in $\ofC'$ because $\oh$ is flat. 
If $V_i$ is a closed point in $\ofC''$, then $U_i$ is a closed point of 
$\ofC'$. We claim that the covering $\ofC'=\cup_i U_i$ 
satisfies property ii) of Proposition \ref{p-3}. 

\begin{lem}\label{l-2}
Let $\beta\in \ofC''$ and $\alpha=\oh(\beta)\in \ofC'$. 
Then $h:\spe^{-1}(\beta)\rightarrow \spe^{-1}(\alpha)$ is surjective.
\end{lem}
{\it Proof.}\ Let $\hat{\co}_\beta$ (resp. $\hat{\co}_\alpha$) 
be the completion 
of the local ring of $\fC''$ at $\beta$ (resp. of $\fC'$ at $\alpha$). 
The natural map $\hat{\co}_\alpha\rightarrow \hat{\co}_\beta$ is finite
because $h$ is quasi-finite.  
Since $\hat{\co}_\alpha$ and $\hat{\co}_\beta$ are integral 
models of $\spe^{-1}(\alpha)$ and $\spe^{-1}(\beta)$ respectively, 
then $h:\spe^{-1}(\beta)\rightarrow \spe^{-1}(\alpha)$ is finite. 
Hence it is surjective because it is flat by 2). 
\hfill $\Box$

Let $x\in \spe^{-1}(U_i)$. By Lemma \ref{l-2}, there exists $y\in 
\spe^{-1}(V_i)\subset \fC''_{K'}$ such that $h(y)=x$. So $Y''_y=Y_x$.  
We deduce that $\# \pi_0(Y_x)$ is constant for $x\in \spe^{-1}(U_i)$. 
The claim and Proposition \ref{p-3} follow.
\hfill $\Box$

\begin{prop}\label{p-2} 
There exists a covering $C=\cup_{i=1}^n C_i$
such that 
\begin{itemize}
\item[{\rm i)}] for each integer $1\leq i\leq n$, either $C_i=C^-(a_{i},
b_i)$ for rational numbers $a_i<b_i$, or $C_i=C^+(a_i,a_i)$
for a rational number $a_i$;
\item[{\rm ii)}] the cardinality of $\pi_0(Y_x)$ is constant 
for $x\in C_i(\Omega)$.
\end{itemize}
\end{prop}
{\it Proof.}\ We start with the covering 
$\cup_{i=1}^n \spe^{-1}(U_i)$ of $C_{K'}$ 
given by Proposition \ref{p-3}, and change it. 
We may assume the $U_i$ connected. We have to consider three cases:
\begin{itemize}
\item[a)] $U_i$ is a closed point $x$ of $\ofC'$ and 
$\spe^{-1}(x)=C^-(s,r)$;
\item[b)] $U_i$ is a closed point $x$ of $\ofC'$ and 
$\spe^{-1}(x)\subset B^-(a,r)\subset C_{K'}$;
\item[c)] $U_i$ is an open subscheme of $\ofC'$.
\end{itemize}
We don't change $U_i$ in case a). 
In case b), since $B^-(a,r)\subset C$, then $v(a)\leq r$. 
Therefore $B^-(a,r)\subset C^+(v(a),v(a))$. Moreover, $\# \pi_0(Y_x)$ 
is constant for $x\in C^+(v(a),v(a))$. So, we replace $\spe^{-1}(U_i)$
by $C^+(v(a),v(a))$. 
In case c), $\spe^{-1}(U_i)$ is a connected affinoid sub-domain 
of $C$.

\begin{lem}\label{l-1}
Let $U$ be a connected affinoid sub-domain of $C$. 
Then, there exist rational numbers 
$0\leq s\leq r \leq e$ such that $U\subset C^{+}(s,r)$ and 
$U$ contains a point of valuation $\alpha$ for any $\alpha
\in {\Bbb Q}\cap [s,r]$. 
\end{lem} 
{\it Proof.}\ It follows from \cite{bgr} 9.7.2/2. 
\hfill $\Box$

Since $\# \pi_0(Y_x)$ depends only on $v(x)$, 
Lemma \ref{l-1} shows that we may replace the affinoid $\spe^{-1}(U_i)$
by $C^+(r,s)$. Proposition \ref{p-2} follows.
\hfill $\Box$

\vspace{2mm}

To finish the proof of Theorem \ref{cont}, we establish the left 
continuity of $\#\pi_0(X^r)$. It may also be possible to deduce 
this result from \cite{blr4}. We give a more direct proof. 
Let $X^{+\infty}$ be the zero locus of $f_1,\dots,f_m$.

\begin{lem}\label{max-p}
{\rm i)}\ Let $g_1,\dots,g_s\in A$ be holomorphic functions, 
$r\in {\Bbb Q}\cup \{+\infty\}$ and $\varepsilon_0$ be 
a rational number.  
If $X(\theta^{\varepsilon_0}g_1,\dots,\theta^{\varepsilon_0}g_s)
\cap X^r=\emptyset$, then there exist rational numbers 
$\varepsilon>\varepsilon_0$ and $a<r$, such that 
$X(\theta^{\varepsilon}g_1,\dots,\theta^{\varepsilon}g_s)\cap X^a
=\emptyset$.\\[3mm]
{\rm ii)}\ Let $U$ be an affinoid sub-domain of $X$ and 
$r\in {\Bbb Q}\cup \{+\infty\}$ such that $U\cap X^r=\emptyset$. 
Then, there exists a rational number $a < r$ such that 
$U\cap X^a=\emptyset$. 
\end{lem}
{\it Proof.}\ i) For $x\in X(\Omega)$, 
we put 
\[
\alpha(x)=\max_{1\leq i\leq n}|f_i(x)|.
\]
By assumption, $\alpha(x)>\theta^r$ for 
$x\in X(\theta^{\varepsilon_0}g_1,\dots,\theta^{\varepsilon_0}g_s)$. 
Since $\alpha$ assumes its minimum on any affinoid variety
(by \cite{bgr} 7.3.4/8), then there is a rational number $r_0<r$ such that 
$\alpha(x)>\theta^{r_0}$ for 
$x\in X(\theta^{\varepsilon_0}g_1,\dots,\theta^{\varepsilon_0}g_s)$. 
For $x\in X(\Omega)$, we put 
\[
\beta(x)=\max_{1\leq i\leq n,1\leq j\leq s}(\theta^{-r_0}|f_i(x)|,
\theta^{\varepsilon_0}|g_j(x)|).
\]
By assumption, $\beta(x)>1$ for all $x\in X(\Omega)$. Thus 
by \cite{bgr} 7.3.4/8, there is a rational number $\delta<0$ such that 
$\beta(x)>\theta^\delta$ for all $x\in X(\Omega)$.
We take $a=r_0+\delta$ and $\varepsilon=\varepsilon_0-\delta$. 
The proof of ii) is similar. \hfill $\Box$

\begin{lem} \label{runge}
Let $U_1,\dots,U_n$ and $V_1,\dots, V_m$ be affinoid sub-domains 
of $X$, and $r\in {\Bbb Q}\cup \{+\infty\}$.
We put $U=\bigcup_{i=1}^nU_i$ and $V=\bigcup_{i=1}^mV_i$,
and we assume that $U\cup V\supset X^r$ and $U\cap V\cap X^r=\emptyset$. 
Then, there exist two admissible open subsets $U'$ and $V'$ of $X$ 
for the strong $G$--topology and a rational number $a < r$, 
such that $U\subset U'$, $V\subset V'$, and 
$U'\cup V'\supset X^a$ and $U'\cap V'\cap X^a=\emptyset$.
\end{lem}
{\it Proof.}\ 
By \cite{bgr} 7.3.5/3, there exists a covering 
$X=\cup_{k=1}^l X_k$ by rational sub-domains $X_k\subset X$ such that 
$U_{i,k}=U_i\cap X_k$ and $V_{j,k}=V_j\cap X_k$ are Weierstrass domains in
$X_k$ for all $1\leq i\leq n$, $1\leq j\leq m$ and $1\leq k\leq l$. 
We put   
$U_{i,k}=X_k(g_{i,k,1},\dots,g_{i,k,s_{i,k}})$ 
and $V_{j,k}=X_k(h_{j,k,1},\dots,h_{j,k,t_{j,k}})$,
and for $\varepsilon\in {\Bbb Q}$, 
\begin{eqnarray*}
U_{i,k}^\varepsilon&=&X_k(\theta^\varepsilon g_{i,k,1},\dots,
\theta^\varepsilon g_{i,k,s_{i,k}}),\\
V_{j,k}^\varepsilon &=&X_k(\theta^\varepsilon h_{j,k,1},\dots,
\theta^\varepsilon h_{j,k,t_{j,k}}).
\end{eqnarray*} 
We put $U^\varepsilon=\bigcup_{i,k} U_{i,k}^\varepsilon$ and 
$V^\varepsilon=\bigcup_{j,k} V_{j,k}^\varepsilon$, 
which are admissible open
subsets of $X$ for the strong $G$--topology by \cite{bgr} 9.1.4/4. 

First, we prove that there exist rational numbers $\varepsilon>0$
and $a_1< r$ such that 
$U^\varepsilon\cap V^\varepsilon\cap X^{a_1}=\emptyset$. 
It is enough to prove that for any $(i,k)$ and $(j,k')$ as above, 
there exist rational numbers $\varepsilon>0$
and $a_1< r$ such that 
\begin{equation}\label{step-1} 
U_{i,k}^\varepsilon\cap V^\varepsilon_{j,k'}\cap X^{a_1}
=\emptyset.
\end{equation}
We work over the affinoid sub-domain $X_k\cap X_{k'}$ of $X$, 
which contains the left hand side of (\ref{step-1}).
Since 
\[
U_{i,k}^\varepsilon\cap V_{j,k'}^\varepsilon=(X_k\cap X_{k'})(
\theta^\varepsilon g_{i,k,1},\dots,\theta^\varepsilon g_{i,k,s_{i,k}},
\theta^\varepsilon h_{j,k',1},\dots,\theta^\varepsilon h_{j,k',t_{j,k'}})
\]
and $U_{i,k}^0\cap V_{j,k'}^0\cap X^r=\emptyset$,
then (\ref{step-1}) follows from Lemma \ref{max-p}-i).

Second, we prove that there exists a rational number $a_1\leq a < r$
such that $U^\varepsilon\cup V^\varepsilon \supset X^a$. 
It is enough to prove that  for any $1\leq k\leq l$, 
there exists a rational number $a_1\leq a < r$
such that 
$(\bigcup_iU_{i,k}^\varepsilon)\cup 
(\bigcup_j V_{j,k}^\varepsilon) \supset X^a\cap X_k$. 
By assumption, $(\bigcup_iU_{i,k}) 
\cup (\bigcup_jV_{j,k})\supset X^r\cap X_k$. 
Therefore, for any choices of 
$g_i\in \{g_{i,k,1},\dots,g_{i,k,s_{i,k}}\}$ for $1\leq i\leq n$,
and $h_j\in \{h_{j,k,1},\dots,h_{j,k,t_{j,k}}\}$ for $1\leq j\leq m$, 
we have 
\[
X_k(\theta^{-\varepsilon}\frac{1}{g_1},\dots,
\theta^{-\varepsilon}\frac{1}{g_n},\theta^{-\varepsilon}\frac{1}{h_1},
\dots,\theta^{-\varepsilon}\frac{1}{h_m})
\cap X^r=\emptyset.
\]
Then by Lemma \ref{max-p}-ii), there exists a rational number 
$a_1\leq a < r$ such that for any choices as above, we have
\[
X_k(\theta^{-\varepsilon}\frac{1}{g_1},\dots,
\theta^{-\varepsilon}\frac{1}{g_n},\theta^{-\varepsilon}\frac{1}{h_1},
\dots,\theta^{-\varepsilon}\frac{1}{h_m})
\cap X^a=\emptyset.
\]
We deduce that $X_k\cap X^{a}\subset (\bigcup_iU_{i,k}^\varepsilon)\cup 
(\bigcup_j V_{j,k}^\varepsilon)$. 
\hfill$\Box$

\vspace{2mm}

Theorem \ref{cont} is a consequence of Proposition \ref{p-2} and the 
following:

\begin{prop}\label{inj} 
For  $r\in {\Bbb Q}\cup \{+\infty\}$, the following natural map
is bijective
\[
\iota_r:\pi_0(X^r)\longrightarrow 
\lim_{\stackrel{\longleftarrow}{a\in {\Bbb Q}_{<r}}}\pi_0(X^a).
\]
\end{prop}
{\it Proof.}\ The injectivity of $\iota_r$
follows from Lemma \ref{runge}.
Observe that a geometrically connected component of an affinoid variety
is an affinoid sub-domain (\cite{bgr} 9.1.4/8 and the discussion after).  
Suppose that $\iota_r$ is not surjective.  
Then, there exists an affinoid sub-domain $U$ of $X$ 
such that $U\cap X^a\not=\emptyset$ for all rational numbers $a<r$,
but $U\cap X^r=\emptyset$. We get a contradiction with Lemma 
\ref{max-p} ii). 
\hfill $\Box$

\section{Ramification of complete intersection rings}\label{lci}

Let $A$ and $B$ be finite flat $\co_K$-algebras, and 
let $u:A\rightarrow B$ be a morphism of $\co_K$-algebras that makes 
$B$ a relative complete intersection over $A$, finite and flat 
of rank $r$ (EGA IV 19.3.6). We fix closed immersions 
$\Spec(A)\rightarrow {\Bbb A}^n_{\co_K}$ and
$\Spec(B)\rightarrow {\Bbb A}^m_{A}$. 
They induce the following commutative diagramme
\[
\xymatrix{
\Spec(B)\ar[r]&{\Bbb A}^m_A\ar[d]\ar[r]&{\Bbb A}^{n+m}_{\co_K}\ar[d]^{\rm pr}\\
& \Spec(A)\ar[r]&{\Bbb A}^n_{\co_K}}
\]
Let $I$ be the ideal sheaf defining the closed immersion $\Spec(B)
\rightarrow {\Bbb A}^m_A$. Since $B$ is semi-local, 
then by EGA IV 19.3.7 there exist an open subscheme $U\subset 
{\Bbb A}^{n+m}_{\co_K}$ and $f_1,\dots,f_m\in \Gamma(U,\co)$ 
such that $\Spec(B)\subset U\cap {\Bbb A}^{m}_A$ and 
\[
\Gamma(U\cap {\Bbb A}^m_A,I)=(f_1,\dots,f_m)\Gamma(U\cap {\Bbb A}^m_A,\co).
\]
Hence, the following diagram is Cartesian
\[
\xymatrix{
\Spec(B)\ar[r]\ar[d]&U\ar[d]^\varphi\\
\Spec(A)\ar[r]& {\Bbb A}^{n+m}_{\co_K}}
\]
where the bottom map is obtained by composing 
the closed immersions $\Spec(A)\rightarrow {\Bbb A}^n_{\co_K}$
and the zero section of ${\rm pr}$, and $\varphi$ is defined 
by $\varphi(x)=({\rm pr}(x),f_1(x),\dots,f_m(x))$. 
Let $\hat{\varphi}_K:\hat{U}_K\rightarrow D^{n+m}$ be the morphism 
induced by $\varphi$ on the rigid fibers of the formal completions 
of $U$ and ${\Bbb A}^{n+m}_{\co_K}$ along their special fibers.
Let $Z$ (resp. $Z'$) be the finite system of generators of $A$
(resp. $B$) over $\co_K$ induced by the closed immersion 
$\Spec(A)\rightarrow {\Bbb A}^n_{\co_K}$ 
(resp. $\Spec(B)\rightarrow {\Bbb A}^{n+m}_{\co_K}$). 
For a rational number $a>0$, we have $X^a_{Z'}\subset \hat{U}_K$
and the following diagram is Cartesian 
\[
\xymatrix{
X^a_{Z'}\ar[r]\ar[d]_{\psi^a}& \hat{U}_K\ar[d]^{\hat{\varphi}_K}\\
X_Z^a\times D^{m,(a)}\ar[r]& D^{n+m}}
\]

\begin{prop}\label{finitude}
For any rational number $a>0$, the morphism 
$\psi^a:X^a_{Z'}\rightarrow X^a_Z\times D^{m,(a)}$ 
is finite and flat of degree $r$. 
\end{prop}
{\it Proof.}\ Let $\oU$ be the special fiber of $U$, and 
let $\spe:D^{n+m}(\Omega)\rightarrow {\Bbb A}^{n+m}(\oOmega)$
and $\spe:\hat{U}_K(\Omega)\rightarrow \oU(\oOmega)$  
be the specialization maps. Let 
\[
T_A=\spe^{-1}({\rm Spec}(A\otimes_{\co_K}\oK))\ 
{\rm and}\ T_B=\spe^{-1}({\rm Spec}(B\otimes_{\co_K}\oK)).
\]
Then $\hat{\varphi}_K$ maps $T_B$ to $T_A$, $X^a_{Z'}\subset T_B$, 
$X^a_{Z}\times D^{m,(a)} \subset T_A$,
and $X^a_{Z'}=\hat{\varphi}_K^{-1}(X^a_Z\times D^{m,(a)})$.
Thus, it is enough to prove 
that $\hat{\varphi}_K$ induces a finite flat 
morphism of rigid analytic spaces 
$T_B\rightarrow T_A$ of degree $r$.

Let $A'$ (resp. $B'$) be the completion of the local rings of 
${\Bbb A}^{n+m}_{\co_K}$ (resp. $U$) at the closed points of
${\rm Spec}(A)$ (resp. ${\rm Spec}(B)$). Then $A'$ and $B'$ provide
integral models of, respectively, $T_A$ and $T_B$ (\cite{dejong} section 7).  
Thus, it is sufficient to prove that the natural morphism 
$A'\rightarrow B'$ makes $B'$ a finite flat $A'$ algebra of rank $r$.
The rings $A'$ and $B'$ are complete regular semi-local Noetherian rings
of dimension $n+m+1$.
Let $\fn$ be the radical of $A'$, $\fm$
be the radical of $B'$, and $(b_1,\dots,b_r)$ be a lifting to $B'$  
of a finite system of generators of the $A$-module $B$. 
By Nakayama, the reductions of $(b_1,\dots,b_r)$ 
generate $B'/\fm^n$ as an $A'/\fn^n$-module. 
Thus, by taking limits, $(b_1,\dots,b_r)$ generate 
the $A'$-module $B'$. 
Since the map $A'\rightarrow B'$ 
is injective, it is finite and flat by EGA $0_{\rm IV}.17.3.5$ (ii). 
\hfill $\Box$

\begin{prop} \label{section} 
Let $A$ be a finite flat $\co_K$--algebra such that $A$ is a complete 
intersection ring (EGA IV 19.3.1) 
and $L=A\otimes_{\co_K} K$ is \'etale over $K$. 
\begin{itemize}
\item[{\rm (i)}] For a rational number $a>0$, the map $\cF(L)\rightarrow 
\cF^a(A)$ is surjective.  
\item[{\rm (ii)}] The map $\displaystyle \cF(L)\rightarrow 
\lim_{\stackrel{\longleftarrow}{a\in {\Bbb Q}_{>0}}}\cF^a(A)$
is bijective. 
\item[{\rm (iii)}] 
Let $E$ be a finite separable extension of $K$ such that $A$ is a 
flat $\co_E$-algebra. Then, the following diagram is cocartesian
\[
\xymatrix{
\cF(L)\ar[r]\ar[d]& \cF^a(A) \ar[d]\\
\cF(E)\ar[r]& \cF^a(\co_E)}
\]
\end{itemize}
\end{prop}
{\it Proof.}\ By EGA IV 19.3.2, $A$ is a 
relative complete intersection over $\co_K$.
(i) follows from Proposition \ref{finitude}. 
(ii) follows from Proposition \ref{inj} for $r=+\infty$.
In order to prove (iii), we apply Proposition \ref{finitude} 
with $A:=\co_E$ and $B:=A$. 
Let $y,y'\in \cF(E)$ with the same image in $\cF^a(\co_E)$. So, there exists 
a geometric connected component $U$ 
of $X^a_{Z}\times D^{m,a}$ that contains both of them. 
By Proposition \ref{finitude},
there exists a geometric connected component $V$ of $X^a_{Z'}$ 
and a finite flat morphism $V\rightarrow U$. 
Let $x,x'\in V(\Omega)$ be inverse images of respectively $y$ and $y'$. 
Then $x,x'\in \cF(L)$, and they have the same image in $\cF^a(A)$. 
The Proposition follows.
\hfill $\Box$

\begin{defi} Let $A$ be a finite flat $\co_K$--algebra
such that $A$ is a complete 
intersection ring and $L=A\otimes_{\co_K} K$ is \'etale over $K$. 
We say that the ramification of $A/\co_K$ is bounded 
by $a$ if the map $\cF(L)\rightarrow \cF^a(A)$ is bijective. 
\end{defi}

\begin{prop}\label{fact-2}
Let $A$ be a finite flat $\co_K$--algebra such that $A$ is a complete 
intersection ring and $L=A\otimes_{\co_K} K$ is \'etale over $K$. 
For a real number $a>0$, we put 
\[
\cF^{a-}(A)=
\lim_{\stackrel{\longleftarrow}{0<b<a}}\cF^b(A)\ \ \ {\rm and} \ \ \ 
\cF^{a+}(A)=
\lim_{\stackrel{\longrightarrow}{b>a}}\cF^b(A),
\]
where $b$ denotes a rational number. Then,  
\begin{itemize}
\item[{\rm (i)}] $\cF^a(A)$ is left-continuous and 
its jumps are rational, i.e., 
$\cF^{a-}(A)=\cF^a(A)$ if $a$ is rational and
$\cF^{a-}(A)=\cF^{a+}(A)$ if $a$ is not rational.
In particular, the number 
\[
c(A/\co_K)=\inf\{ a\in {\Bbb Q}_{>0}\ ; \cF(L)\simeq \cF^a(A)\}
\]
is rational. We call it the conductor of the extension $A/\co_K$;
\item[{\rm (ii)}] The ramification of $A/\co_K$ is not bounded by 
$c(A/\co_K)$;
\item[{\rm (iii)}] The extension $A/\co_K$ is \'etale if and only if 
$c(A/\co_K)=0$. 
\end{itemize}
\end{prop}
{\it Proof.}\ (i) follows from Theorem \ref{cont} and Proposition 
\ref{section}-(i).
(ii) is a consequence of the left-continuity of $\cF^a(A)$. 
(iii) follows from Lemma \ref{extra}-ii), since $A/\co_K$ 
is \'etale if and only if the map 
$\cF(L)\rightarrow \Hom_{\co_K}(A,\oOmega)$ is bijective.
\hfill $\Box$

\begin{lem}\label{b-c} Let $K'$ be a finite separable extension of $K$ 
of ramification index $e$, contained in $\Omega$. 
Let $A$ be a finite flat $\co_{K}$--algebra, 
$B$ be a finite flat $\co_{K'}$--algebra, and 
$u:A\rightarrow B$ be a morphism of $\co_K$-algebras. For a 
rational number $a>0$, we denote by $\cF_{K'}^a$ 
the functor constructed in Subsection $\ref{construction-1}$ 
for finite flat $\co_{K'}$-algebras. Then, we have a natural map 
\[
\cF^{ea}_{K'}(B)\rightarrow \cF^a(A),
\]
that is an isomorphism if $B=A\otimes_{\co_{K}}\co_{K'}$.\\[2mm]
Assume furthermore that $A$ and $B$ are complete intersection 
rings, that $L=A\otimes_{\co_K}K$ 
is \'etale over $K$ and $L'=B\otimes_{\co_{K'}}K'$ 
is \'etale over $K'$, and that the map 
$\Hom_{K'}(L',\Omega)\rightarrow \Hom_{K}(L,\Omega)$ induced by $u$ 
is injective. If the ramification of $A/\co_{K}$ is
bounded by $a$, then the ramification of $B/\co_{K'}$ is
bounded by $ea$
\end{lem}
{\it Proof.}\ Let $Z$ and $Z'$ be finite systems of generators
of, respectively, $A/\co_K$ and $B/\co_{K'}$, such that $u(Z)\subset Z'$. 
Let $X^a_Z$ and $Y^{ea}_{Z'}$
be the affinoid varieties associated with, respectively, $(A/\co_K,Z)$ 
and $(B/\co_{K'},Z')$. Then, the natural projection 
induces a rigid-analytic morphism 
$Y^{ea}_{Z'}\rightarrow X^a_{Z}\times_{K}K'$. 
We deduce a morphism $\cF^{ea}_{K'}(B)\rightarrow \cF^a(A)$. 
If $B=A\otimes_{\co_K}\co_{K'}$, we take $Z'=Z\otimes 1$. 
Hence, $Y^{ea}_{Z'}\simeq X^a_{Z}\times_{K}K'$. 

We have a natural commutative diagram
\[
\xymatrix{
\Hom_{K'}(L',\Omega)\ar[r]^i\ar[d]_{g}&\cF^{ea}_{K'}(L')\ar[d]\\
\Hom_K(L,\Omega)\ar[r]_j&\cF^{a}(L)}
\]
in which $g$ is injective, and $i$ is surjective. Since 
$j$ is bijective, we deduce that $i$ is also bijective. 
\hfill $\Box$

\subsection{Example~: monogenic extensions}\label{monogenic}

Let $L$ be a finite separable extension of $K$ of degree $d$.  
Assume its valuation ring $\co_L$ monogenic over 
$\co_K$, and fix $\co_L=\co_K[X]/P$ a monogenic presentation
(so $P$ is a monic polynomial of degree $d$). Let 
$z_1,\dots,z_d$ be the zeros of $P$ in $\co$. 

\begin{lem}\label{newton}
Let $a>0$ be a rational number. The following are equivalent~:
\begin{itemize}
\item[{\rm (i)}] the ramification of $\co_L/\co_K$ is bounded by $a$;
\item[{\rm (ii)}] $\displaystyle \sum_{i\not=1}v(z_i-z_1)
+\sup_{i\not=1}v(z_i-z_1)<a$.
\end{itemize}
\end{lem}
{\it Proof.}\ Notice that the left hand side of inequality 
(ii) does not depend on the numbering of the zeros of $P$
because ${\rm Gal}(\Omega/K)$ acts 
transitively on them. 
First, assume that (ii) does not hold. Then 
$D(z_1,\theta^\alpha) \subset P^{-1}(D(0,\theta^a))$,
where $\alpha=\sup_{i\not=1}v(z_i-z_1)$. 
But $D(z_1, \theta^\alpha)$ is connected and contains 
at least two zeros of $P$. Therefore, 
the ramification of $\co_L/\co_K$ is not bounded by $a$.
Second, assume that (ii) holds. Let $y\in \co$ be
such that $v(P(y))\geq a$. 
We may assume that $v(y-z_1)\geq v(y-z_i)$
for $2\leq i\leq d$. Then, 
\[
v(P(y))\leq v(y-z_1)+\sum_{i\not=1}v(z_i-z_1). 
\]
We deduce that $v(y-z_1)\geq \alpha+\varepsilon$, where $\varepsilon>0$ 
is the difference between the right and the 
left hand sides of inequality (ii). 
Hence, $P^{-1}(D(0,\theta^a))=\bigsqcup_{1\leq i\leq d}D(z_i,\theta^{\alpha
+\varepsilon})$, and the ramification of $\co_L/\co_K$ is bounded by $a$. 
\hfill $\Box$

\vspace{0.5cm}

An important example of monogenic valuation ring extensions is provided 
by unfiercely ramified extensions, i.e. finite separable extensions 
of local fields with separable residue extensions
(see \cite{serre1} III 6 Proposition 12). 
For these extensions, our theory gives the classical 
ramification theory developed in \cite{serre1,deligne}, which we 
summarize below. Let $L/K$ be 
a finite separable unfiercely ramified extension, 
$\cF(L)={\rm Hom}_K(L,\Omega)$, and $u\geq -1$ be a real number. We define 
an equivalence relation $R_u$ on $\cF(L)$ by 
\[
\sigma\equiv \tau \ ({\rm mod}\ R_u) \ \ \Leftrightarrow 
\ \ v_L(\sigma-\tau):=\inf_{x\in \co_L}v_L(\sigma(x)-\tau(x))\geq u+1,
\]
where $v_L=r v $ and $r$ is the ramification index of $L/K$.
The Galois group ${\rm Gal}(\Omega/K)$ acts transitively on $\cF(L)$ 
and preserves 
the equivalence relation $R_u$. Hence, the equivalence classes have the 
same cardinality $r_u$. The Herbrand function $\varphi_{L/K}:
[-1,+\infty)\rightarrow[-1,+\infty)$ is defined by
\[
\varphi_{L/K}(u)=\int_0^u\frac{r_t}{r_0}dt.
\]
Let $\psi_{L/K}$ be 
its inverse. The upper numbering 
equivalence relation on $\cF(L)$ is defined by $R^f=R_{\psi_{L/K}(f)}$.

\begin{prop}\label{classical}
Let $f>-1$ be a rational number. 
The following are equivalent~:
\begin{itemize}
\item[{\rm (i)}] the ramification of $\co_L/\co_K$ is bounded by $f+1$;
\item[{\rm (ii)}] the equivalence relation $R^f$ is trivial; 
\end{itemize}
\end{prop}
{\it Proof.}\ Let $\sigma_0\in \cF(L)$. For $\sigma\not=\sigma_0$, 
put $v(\sigma-\sigma_0)=\frac 1 r v_L(\sigma-\sigma_0)$. 
By \cite{deligne} Proposition A.6.1, (ii) is equivalent
to 
\[
\sum_{\sigma\not=\sigma_0}v(\sigma-\sigma_0)+
\sup_{\sigma\not=\sigma_0}v(\sigma-\sigma_0)<f+1,
\]
which is equivalent to (i) by Lemma \ref{newton}. 
\hfill $\Box$

\subsection{Proofs of \ref{fil1}, \ref{properties} and \ref{rational}}
\label{ram-fil1}

Theorem \ref{fil1} follows from Proposition 
\ref{section}. 
Theorem \ref{rational} follows from Proposition \ref{fact-2}-(i). 
Proposition \ref{properties}-3) was proved in Subsection \ref{monogenic}. 
We prove Proposition \ref{properties}-2).
The inclusion $G_{K'}^{ma}\subset G^a$ 
is a consequence of Lemma \ref{b-c}. 
Assume that $K'/K$ is unramified, and let $L$ be a finite separable extension
of $K$. Then by Lemma \ref{b-c}, we have $\cF^{a}_{K'}(L\otimes_K K')=
\cF^a(L)$, so $\cF(L)/G^a=\cF(L)/G^{a}_{K'}$. 
We deduce that $G^a_{K'}=G^a$.
Proposition \ref{properties}-1) is a consequence of the following 
two Propositions.

\begin{prop}\label{tame}
Let $L$ be a finite separable extension of $K$. 
The following are equivalent~:
\begin{itemize}
\item[$1)$] $L/K$ is tamely ramified;
\item[$2)$] the ramification of $\co_L/\co_K$ is bounded by $a$, 
for any rational number  $a>1$.
\end{itemize}
\end{prop}
{\it Proof.} \ $1)\Rightarrow 2)$.\ Follows from classical ramification 
theory. 
$2) \Rightarrow 1)$.\  
We assume first that $\co_L/\co_K$ is monogenic. 
Lemma \ref{newton} implies that   
\[
\sum_{i\not=1}v(z_i-z_1)+
\sup_{i\not=1}v(z_i-z_1)\leq 1
\]
If $\sup_{i\not=1}v(z_i-z_1)=0$ then $\co_L/\co_K$ 
is \'etale. Otherwise, $v({\frak d}_L)=
\sum_{i\not=1}v(z_i-z_1)<1$ where ${\frak d}_L$ is the 
different of $\co_L/\co_K$. We conclude using Proposition \ref{car-tame}. 

The proof in the general case is by induction on the lexicographical order 
of $(v(\delta_L),s_L)$ (cf. Appendix for the notations). 
The result is trivial if $v(\delta_L)=0$, and proved above if 
$s_L=0$. Assume $s_L\not=0$. By Lemma \ref{b-c}, we may assume 
the residue extension $\oL/\oK$ purely inseparable. 
Let $K'$ be as in  lemma 
\ref{raynaud1}, and  $L'$ be a composed extension of $L$ and $K'$. 
By Lemma \ref{b-c}, assumption 2) implies that the ramification of 
$\co_{L'}/\co_{K'}$ is bounded by $a$ for any $a>1$. 
Since $(v(\delta_{L'}),s_{L'})<
(v(\delta_L),s_L)$ (cf. lemma \ref{raynaud1}), we deduce that 
$L'/K'$ is tamely ramified. Let $e$ be the ramification index of 
$L/K$, $e'$ be the ramification index of $L'/K'$, and $r$ 
be the ramification index of $L'/L$, so $e'=er$. 
We have $\oK\subset \oK'\subset 
\oL\subset \oL'$, moreover $\oL'/\oK'$ is separable and $\oL/\oK$ 
is purely inseparable. It follows that $\oL=\oK'$. 
On the one hand $[L':L]=r [\oL':\oL]$ is a divisor of 
$[K':K]$ which is a power of $p$. On the other hand, 
$r$ is prime to $p$ because so is $e'$. Hence $r=1$. 
We conclude that $\co_{L'}/\co_{L}$ is \'etale. 

Let $F$ be the maximal sub-extension of $L'/K$ such that 
$\co_{F}/\co_K$ is \'etale. Since $\oL/\oK$ is purely inseparable, 
then $L$ and $F$ are disjoint. Moreover, 
$[F:K]$ is the separable factor of $[\oL':\oK]$ which is 
$[\oL':\oL]=[L':L]$. Hence, $L'$ is the composed extension of $L$ and $F$. 
By Proposition \ref{fact-2} (iii), the ramification of $\co_{F}/\co_K$
is bounded by $1$. 
We deduce that the ramification of $\co_{L'}/\co_K$ 
is bounded by $a$, for any $a>1$. Therefore, 
the ramification of $\co_{K'}/\co_K$ 
is bounded by $a$, for any $a>1$.
But $\co_{K'}/\co_K$ is monogenic. Then $K'/K$ is tamely ramified, 
which is a contradiction with $ii)$ of Lemma \ref{raynaud1}. 
\hfill $\Box$

\begin{prop} \label{etale2}
Let $L$ be a finite separable extension of $K$. The following are equivalent~:
\begin{itemize}
\item[$1)$] $\co_L/\co_K$ is \'etale;
\item[$2)$] the ramification 
of $\co_L/\co_K$ is bounded by $a$, for any rational number $a>0$;
\item[$3)$] the ramification 
of $\co_L/\co_K$ is bounded by $1$. 
\end{itemize}
\end{prop}
{\it Proof.}\ The equivalence between 1) and 2) was proved in
Proposition \ref{fact-2}-(iii).  We prove that 3) implies 1). 
Assume that condition 3) holds true. Proposition \ref{tame}
implies that $L/K$ is tamely ramified and in particular 
the extension $\co_L/\co_K$ is monogenic. 
Let $e$ be the ramification index of $L/K$ and ${\frak d}_L$ be a different 
of $\co_L/\co_K$. Lemma \ref{newton} implies that
\[
v({\frak d}_L)+\sup_{i\not=1}v(z_i-z_1)=(1-\frac 1 e) 
+\sup_{i\not=1} v(z_i-z_1)<1. 
\]
Thus $\sup_{i\not=1}v(z_i-z_1)<1/e$. 
If $e\not=1$, then 
\[
v({\frak d}_L)\leq (e-1) \sup_{i\not=1}v(z_i-z_1)<\frac{e-1}{e},
\]      
which is impossible. We deduce that $\co_L/\co_K$ is \'etale. 
\hfill $\Box$

\section{Deformation to the special fiber}

\begin{lem}\label{ci}
Let $L$ be a finite separable extension of $K$. Then 
$\co_L$ is a complete intersection over $\co_K$, i.e. 
$\co_L\simeq \co_K[X_1,\dots,X_n]/(f_1,\dots,f_n)$. 
\end{lem}
{\it Proof.}\ Let $\varphi:R=\co_K[X_1,\dots,X_n]\rightarrow \co_L$
be a surjective morphism of $\co_K$-algebras, $I={\rm ker}(\varphi)$
and $\fm=\varphi^{-1}(\fm_L)$. Let $\hat{R}$ and $\hat{I}$
be the $\fm$-completions of $R$ and $I$.
We have $\co_L=\hat{R}/\hat{I}$. Then by EGA IV 19.3.2, there 
exist $(g_1,\dots,g_n)\in \hat{R}$ such that $\hat{I}=(g_1,\dots,g_n)
\hat{R}$. Thus $I/\fm I$ is generated as an $\co_K$-module 
by the residue classes of $n$ equations $(f_1,\dots,f_n)\in I$.
By Nakayama, we have $\hat{I}=(f_1,\dots,f_n)\hat{R}$. 
We put $D=R/(f_1,\dots,f_n)$ and $H=I/(f_1,\dots,f_n)$, so we have the exact 
sequence $0\rightarrow H\rightarrow D\rightarrow \co_L\rightarrow 0$.
By completion at $\fm$, we get that $\hat{H}=0$. Hence 
\[
H/H^2=H\otimes_DD/H=H\otimes_D\hat{D}=\hat{H}=0.
\]
We deduce that $H=H^2$ and $(1-e)H=0$ for some $e\in H$. Then $H=eD$,
$e=e^2$ and 
\[
\co_L\simeq D/H\simeq D[X]/((1-e)X-e)\simeq \co_K[X_1,\dots,X_{n+1}]/
(f_1,\dots,f_{n+1}),
\]
where $f_{n+1}=(1-\tilde{e})X_{n+1}-\tilde{e}$ and 
$\tilde{e}\in R$ is a lifting of $e$. \hfill $\Box$

Let $L$ be a finite separable extension of $K$. 
We fix a presentation of complete intersection
\[
\co_L=\co_K[X_1,\dots,X_n]/(f_1,\dots,f_n),
\] 
and denote by $Z=(z_1,\dots,z_n)$ the residue classes of $(X_1,\dots,X_n)$.
Let $f:D^n\rightarrow D^n$ 
be the map defined by the equations $(f_1,\dots,f_n)$, 
and for a rational number $a>0$,
let $f^a:X^a_Z\rightarrow D^{n,(a)}$ be the induced map. 
Let $\fX^a_Z$ (resp. $\cD^{n,(a)}$) be the normalized integral 
model of $X^a_Z$ (resp. $D^{n,(a)}$) over $\co$, and let 
$\ofX^a_Z$ (resp. $\ocD^{n,(a)}$) be its special fiber. 
By the finiteness theorem of Grauert and Remmert (\cite{gr} and \cite{blr4} 
theorem 1.2), the map $f^a$ extends to a finite morphism 
$f^a:\fX^a_Z\rightarrow \cD^{n,(a)}$. We denote by $\of^a:
\ofX^a_Z\rightarrow \ocD^{n,(a)}$ the induced morphism
on the special fibers. 
In the sequel, we fix $Z$ and drop it from the notation.
The goal of this section is to prove the following

\begin{teo}\label{t-def}
Assume that $L/K$ is not tamely ramified, and let $c=c(L/K)>1$.
Then $f^c:\fX^c\rightarrow \cD^{n,(c)}$ is \'etale.
\end{teo}

\begin{prop}\label{dif-ram}
Let $e$ be the ramification index of $L/K$, 
and let $\alpha_1,\dots,\alpha_r$ be the positive integers such that
$\Omega^1_{\co_L/\co_K}\simeq \oplus_{i=1}^r \co_{L}/\pi_L^{\alpha_i}
\co_L$. If $\co_L/\co_K$ is not \'etale, 
then $\alpha_i<e c(L/K)$ for all $i$. 
\end{prop}
{\it Proof.}\ We may assume the residue extension $\oL/\oK$ 
purely inseparable. We fix a presentation as in Lemma \ref{ci}, and denote
by $f'\in {\rm M}_n(\co_K[X_1,\dots,X_n])$ 
the Jacobian matrix of the equations $(f_1,\dots,f_n)$ 
and $J\in {\rm M}_n(\co_L)$ its residue class. 
Let $\sigma:\co_L\rightarrow \co$ be an embedding.
To prove the Proposition, it is enough to prove that there 
exists a rational number $\gamma>0$ 
such that if the ramification of $\co_L/\co_K$ is bounded by $a>0$,
then $a>\gamma$ and $\pi^{a-\gamma}\co^n \subset \sigma(J)\co^n$. 
We need the following:
 
\begin{lem}\label{approx}
For any rational number $a>0$, there exists a rational number 
$\alpha>0$, such that $X^a\subset \bigcup_{\sigma\in \cF(L)}
(\sigma+D^{n,(\alpha)})$.
\end{lem}
{\it Proof.}\ Let $0<b<a$ be a rational number, and consider the following 
commutative diagram
\[
\xymatrix{
X^a\ar[r]\ar[d]& X^b\ar[r]^{\spe}\ar[d]^{f^b}&
{\ofX^b}\ar[d]^{\of^b}\\
D^{n,(a)}\ar[r]& D^{n,(b)}\ar[r]^{\spe}& {\ocD}^{n,(b)}}
\]
Since $\of^b$ is finite and $\spe(D^{n,(a)})$ is a closed point, 
then $\spe(X^a)$ is a finite set of closed points. 
So, if we consider $X^a$ as an affinoid sub-domain 
of $D^n$ and denote by $\spe:D^n\rightarrow \ocD^n$
the specialization map, then $\spe(X^a)$ is also a finite set of closed 
points. The Lemma follows using the maximum modulus principle.
\hfill $\Box$

The assumption that $\oL/\oK$ is purely inseparable and Lemma \ref{approx} 
imply that there exists a rational number $0<\gamma<1$
such that $X^1\subset \sigma+D^{n,(\gamma)}$. 
Let $a>0$ be a rational number such that the ramification 
of $\co_L/\co_K$ is bounded by $a$.
Observe that $a>1$ because $\co_L/\co_K$ is not \'etale (by  
Proposition \ref{etale2}). 
Let $s:D^{n,(a)}\rightarrow X^a$ be the section of $f^a$ that maps 
$0$ to $\sigma$, and $\psi$ be 
the composite of $s$ and the canonical embedding $X^a
\subset D^n$. We write $\psi=(g_1,\dots,g_n)$, where 
$g_i=\sum_\nu a_{i,\nu} X^\nu\in \Omega [[X_1,\dots,X_n]]$ and satisfies
$\lim_{\nu\mapsto +\infty} v(a_{i,\nu})+a|\nu|=+\infty$. 
The inclusion $X^1\subset \sigma+D^{n,(\gamma)}$ implies that
\[
|g_i-g_i(0)|_{\sup}=\sup_{|\nu|\geq 1}|a_{i,\nu}|\theta^{a|\nu|}\leq \theta^\gamma.
\]
Therefore $\psi'(0)(\pi^{a} \co^n) \subset \pi^\gamma \co^n$. 
Proposition \ref{dif-ram} follows by observing that $\psi'(0)$ is the inverse
of $\sigma(J)$. 

\hfill $\Box$

\begin{prop}\label{r-etale}
Let $e$ and $\alpha_1,\dots,\alpha_r$ be as in Proposition $\ref{dif-ram}$,
and let $a>1$ be a rational number such that $ea>\alpha_i$ for all $i$. 
Then $f^{a}:X^{a}\rightarrow D^{n,(a)}$ is \'etale.
\end{prop}
{\it Proof.}\ Let 
$f'\in {\rm M}_n(\co_K[X_1,\dots,X_n])$ be the Jacobian matrix 
of the equations $(f_1,\dots,f_n)$, 
let $x\in X^{a}(\Omega)$ be a zero of $\det(f')$, 
and let $\varphi:\co_L\rightarrow \co/\pi^{a}\co$ be the associated 
morphism. Let $\pi_L$ be a uniformizer of $L$ 
and $v\in \co_L^*$ be a unit such that $\pi_L^e=v\pi$. 
Since $\varphi(\pi_L)^e=\varphi(v) \pi \ ({\rm mod} \ \pi^{a})$ and $a>1$,
then there exists a unit $u\in \co^*$ such that 
$\varphi(\pi_L)=u\pi_L \ ({\rm mod}\ \pi^a)$.  
Since $\alpha_i<ea$ for all $i$, then 
\[
\Omega^1_{\co_L/\co_K}\otimes_{\varphi}\co/\pi^a\co \simeq 
\bigoplus_{i=1}^r\co/\pi_L^{\alpha_i}\co.
\]
On the other hand, we have 
\[
\Omega^1_{\co_L/\co_K}\otimes_{\varphi}\co/\pi^a\co \simeq 
{\rm Coker}(\xymatrix{\co^n\ar[r]^{f'(x)}& 
\co^n})\otimes \co/\pi^a\co.
\]
Since the rank of the matrix $f'(x)$ is $\leq n-1$, then 
$\Omega^1_{\co_L/\co_K}\otimes_{\varphi}\co/\pi^a\co$ has a direct summand 
isomorphic to $\co/\pi^a\co$. 
We get a contradiction with the assumption that $ea>\alpha_i$ for all $i$. 
We deduce that $\det(f'(x))\not=0$ for all $x\in X^{a}(\Omega)$. 
\hfill $\Box$

\vspace{2mm}

{\it Proof of Theorem} \ref{t-def}.\ The map $f^c:X^c\rightarrow D^{n,(c)}$
is \'etale by Propositions \ref{dif-ram} and \ref{r-etale}.
Let $\spe:X^c \rightarrow \ofX^c$ be the specialization
map,  $\overline{0}$ be the center of $\ocD^{n,(c)}$, 
$X^{c+}=\bigcup_{b>c}X^b$, and $D^{n,(c+)}=\bigcup_{b>c}D^{n,(b)}$. 
Then $(f^c)^{-1}(\overline{0})=\spe(\cF(L))\simeq
\pi_0(X^{c+})$, and by definition of $c$,  
$\pi_0(X^{c+})\simeq \cF(L)$ and the map 
$f^{c+}:X^{c+}\rightarrow D^{n,(c+)}$ is finite \'etale and geometrically 
totally decomposed. 
Thus for any $x\in \spe(\cF(L))$, the completion of the local ring 
of $\fX^c$ at $x$ is isomorphic to the completion 
of the local ring of $\cD^{n,(c)}$ at $\overline{0}$
(see \cite{bosch} Section 6, or \cite{bl} Lemma 2.1). 
We deduce that $f^c:\fX^c\rightarrow \cD^{n,(c)}$ is \'etale above
$\overline{0}$, and hence 
above the generic point of the closed fiber of $\cD^{n,(c)}$. 
Since $\fX^c$ is normal and $\cD^{n,(c)}$ is regular,
then by Zariski's purity theorem
$f^c:\fX^c\rightarrow \cD^{n,(c)}$ is \'etale. 
\hfill $\Box$

\vspace{2mm}

The geometric monodromy induces an $\oOmega$-linear action of the inertia group
$G^1$ on $\ofX^c$. 

\begin{cor}\label{d-c+}
The group $G^c$ stabilizes the (geometric) 
connected components of $\ofX^c$, and the subgroup $G^{c+}$ acts trivially 
on $\ofX^c$.
\end{cor}
{\it Proof.}\ We know that $G^c$ stabilizes the geometric connected 
components of $X^c$, and hence those of $\ofX^c$. 
By Lemma \ref{d-mon} below, $G^{c+}$ acts trivially 
on $\ocD^{n,(c)}$ because $c>1$.
Let $\spe:X^c \rightarrow \ofX^c$ be the specialization
map.
Since $\spe(\cF(L))\simeq \pi_0(X^{c+})=\cF(L)$, then 
$G^{c+}$ acts trivially on $\spe(\cF(L))$.
Moreover the map $\ofX^c\rightarrow \ocD^{n,(c)}$ is \'etale, 
and each connected component of $\ofX^c$ contains a point of $\spe(\cF(L))$.
Then $G^{c+}$ acts trivially on $\ofX^c$.
\hfill $\Box$

\begin{lem}\label{d-mon}
For any rational number $a\geq 0$, 
the geometric monodromy induces a trivial action of $G^{1+}$ on $\ocD^{n,(a)}$.
\end{lem}
{\it Proof.}\ Let $E$ be a finite Galois extension of $K$ containing 
an element of valuation $a$. We have $D^{n,(a)}_{E}=\Sp(A_{E})$, where 
\[
A_{E}=\{\sum_\nu a_\nu X^\nu \in E[[X_1,\dots,X_n]]; 
\ \lim_{|\nu|\rightarrow +\infty} |a_\nu|\theta^{a|\nu|}=0\}.
\]
The normalized integral model 
$\cD^{n,(a)}$ of $D^{n,(a)}$ is defined over $\co_{E}$, and we have
$\cD^{n,(a)}_{\co_{E}}=\Spf(A')$ where 
\[
A'=\{\sum_\nu a_\nu X^\nu \in A_{E}; 
\ \ \sup_\nu |a_\nu|\theta^{a|\nu|}\leq 1\}.
\] 
Let $\fm_\Omega$ be the maximal ideal of $\co$. 
Then $G^{1+}=P$ is the kernel of the map
$G \to {\rm Aut}(\Omega^*/1+m_\Omega)$.
Therefore, for $\sigma\in G^{1+}$ and 
$a\in E-\{0\}$, we have $|\sigma(a)-a|<|a|$. The Lemma follows.
\hfill $\Box$

\vspace{2mm}

Suppose further that $L$ is a finite Galois extension of $K$ contained 
in $\Omega$. Let $G_L={\rm Gal}(\Omega/L)$, $G(L/K)=G/G_L$
be the Galois group of $L/K$, and  $(G(L/K)^a, 
a\in {\Bbb Q}_{\geq 0})$ be the quotient filtration of $(G^a, 
a\in {\Bbb Q}_{\geq 0})$. 

\begin{cor}
Let $H$ be a geometric connected component of $\ofX^c$ and let 
$H\rightarrow \ocD^{n,(c)}$ be the restriction of $\of^c$ to $H$. 
Then $H\rightarrow \ocD^{n,(c)}$ is an \'etale Galois covering 
of group $G(L/K)^c$. 
\end{cor}
{\it Proof.}\ We have $G(L/K)^c=G^c/G^c\cap G_L$.
Let $\spe:X^c \rightarrow \ofX^c$ be the specialization
map. The group $G^c\cap G_L$ acts trivially on $\ocD^{n,(c)}$
(by Lemma \ref{d-mon}), and on $\spe(\cF(L))\simeq \cF(L)$. 
As in the proof of Corollary 
\ref{d-c+}, we deduce that $G^c\cap G_L$ acts trivially on $\ofX^c$. 
Thus the action of $G^c$ on $\ofX^c$ factors through 
$G(L/K)^c$. Moreover, $G(L/K)^c$ acts transitively and freely 
on the fiber of $H\rightarrow \ocD^{n,(c)}$ above the origin. 
The Corollary follows. \hfill $\Box$

\section{Ramification and Newton polygons}\label{rnp}

Let $f(X)=a_0X^d+a_1 X^{d-1}+\dots+a_d\in \co_K[X]$ be a polynomial. 
Its Newton polygon $N_f:[0,d]\rightarrow {\Bbb R}\cup \{+\infty\}$
is the convex envelop of the set of points $(i,v(a_i))$ for $0\leq i\leq d$.
We say that $f$ has a Newton polygon of type $0=d_0\leq d_1\leq 
\dots\leq d_r=d$ if the slope of $N_f$ is constant 
on $]d_i,d_{i+1}[$. The type of a Newton polygon
can be described by a finite number of divisibility 
relations between monomials in the coefficients $a_0,\dots,a_d$.
We denote these relations by $D(d_0,\dots,d_r;a_0,\dots,a_d)$.  

\begin{defi}
Let $A$ be a commutative ring and $f(X)=a_0X^d+a_1 X^{d-1}+\dots+a_d\in 
A[X]$ be a polynomial. We say that $f$ has a Newton polygon of type $[d_0,
\dots,d_r]$ if the divisibility relations $D(d_0,\dots,d_r;a_0,\dots,a_d)$
are satisfied. 
\end{defi}

Let $K\subset L\subset M$ be finite separable extensions such that 
\begin{itemize}
\item[(a)] $M/K$ is not tamely ramified, so $c=c(M/K)>1$;
\item[(b)] $\co_M/\co_L$ is not \'etale and has degree $[M:L]=p$;
\item[(c)] $\cF^c(M)\simeq \cF^c(L)$.
\end{itemize}
We fix the presentations (by Lemma \ref{ci})
\begin{eqnarray*}
\co_L&=&\co_K[y_1,\dots,y_n]/(f_1,\dots,f_n):=
\co_K[y]/(f),\\
\co_M&=&\co_L[x]/q(x)=\co_K[y,x]/(f_1,\dots,f_n,g(y,x)),
\end{eqnarray*}
where $y:=(y_1,\dots,y_n)$ and  
$g(y,x)$ is a lifting of $q(x)$ that is monogenic of degree 
$p$ in $x$. We set 
\[
\psi: 
\begin{array}[t]{clcr}
D^n\times D&\longrightarrow & D^n\times D\\
(y,x)&\longmapsto& (y,g(y,x)).
\end{array}
\]
For a rational number $a>0$, let $X^a\subset D^n\times D$
and $Y^a\subset D^n$ be the affinoid varieties associated respectively with 
$M/K$ and $L/K$ and the above presentations.  
We have $X^a=\psi^{-1}(Y^a\times D^{(a)})$. We denote by 
$\psi^a:X^a \rightarrow Y^a\times D^{(a)}$ the restriction of $\psi$.

Let $\fX^a$, $\fY^a$ and $\cD^{(a)}$ be the normalized integral models 
of respectively $X^a$, $Y^a$ and $D^{(a)}$ over $\co$.
By the finiteness theorem of Grauert and Remmert, the map $\psi^a$
extends to a finite map $\psi^a:\fX^a\rightarrow \fY^a\times \cD^{(a)}$. 

\begin{prop}\label{rec1}
The map $\psi^c:\fX^c\rightarrow \fY^c\times \cD^{(c)}$ is \'etale.
\end{prop}
{\it Proof.}\ Indeed, the maps $\fX^c\rightarrow \cD^{n+1,(c)}$ 
and $\fY^c\rightarrow \cD^{n,(c)}$ are \'etale by 
Theorem \ref{t-def}. 
\hfill $\Box$

\vspace{2mm}

Let $Z^a=Y^a\times D^{(a)}$ and let $\cZ^a=\fY^a\times \cD^{(a)}$ 
be its normalized integral model over $\co$.  
We denote by $z$ the canonical coordinate on $D$, and also 
the induced function on $\cZ^a$. We put $P(x)=g(y,x)\in\co_{Y^a}[x]$.
Then we have 
\[
X^a=\Sp\left(\co_{Z^a}[x]/(P(x)-z)\right).
\]
Let $\tau$ be the residue class of $x$, which is a function in $\co_{\fX^a}$.

\begin{teo}\label{main-np-t}
The polynomial $P(x+\tau)-z\in \co_{\fX^c}[x]$ has a Newton polygon 
of type $[0, p-1,p]$.
Moreover, $P'(\tau)\co_{\fX^c}=\beta^{p-1}\co_{\fX^c}$ 
for an element $\beta\in \co$ with $v(\beta)=b/p$,
where $e_{L/K}b=c(M/L)$ and $e_{L/K}$ is the ramification 
index of $L/K$.  
\end{teo}
{\it Proof.}\ Let $Y^c_\circ$ be a geometrically connected component 
of $Y^c$, $Z^c_\circ=Y^c_\circ\times D^{(b)}$ , 
and $X^c_\circ$ be the geometrically connected component 
of $X^c$ above it. Let $K'$ be a finite separable extension of $K$
over which $Y^c_\circ$, $Z^c_\circ$ and $X^c_\circ$ are defined 
as well as their normalized integral models $\fY^c_\circ$, $\cZ^c_\circ$ 
and $\fX^c_\circ$. 
By Theorem \ref{t-def}, the closed fibers $\ocZ^c_\circ$ 
of $\cZ^c_\circ$ and $\ofX^c_\circ$ of $\fX^c_\circ$
are smooth. So they are geometrically reduced and irreducible because 
they are geometrically connected. 
Since $\psi^c:\fX^c_\circ\rightarrow \cZ^c_\circ$ finite flat of degree $p$,
then by Proposition \ref{3.5}   
there exist a function $\xi$ on $\fZ^c_\circ$ 
and $\beta\in \co_{K'}$ such that $(\tau+\xi)/\beta\in
\co_{\fX^c_\circ}$, and   $\fX^c_\circ\simeq 
\Spf(\co_{\fZ^c_\circ}[(\tau+\xi)/\beta])$.

Let $\gamma=(\tau+\xi)/\beta$ and let $F(X)\in \co_{\fZ^c_\circ}[X]$
be its characteristic polynomial. 
We proved in Proposition \ref{3.5} that 
$\fX^c_\circ=\Spf(\co_{\fZ^c_\circ}[X]/F)$. 
Since $P(X)-z$ is the characteristic polynomial of $\tau$, 
then $\beta^p F(X)=P(\beta X-\xi)-z$.
We deduce that $\beta^{p-1} F'(\gamma)=P'(\tau)$ and 
\[
\beta^pF(X+\gamma)=P(\beta X+\tau)-z.
\]
Since $\psi^c:\fX^c_\circ\rightarrow \fZ^c_\circ$ is \'etale
by Proposition \ref{rec1}, then 
$F'(\gamma)\co_{\fX^c_\circ}=\co_{\fX^c_\circ}$
and $F(X+\gamma)$ has a Newton polygon of type $[0,p-1,p]$. 
Therefore, $\beta^{p-1}\co_{\fX^c_\circ}
=P'(\tau)\co_{\fX^c_\circ}$ and $P(X+\tau)-z$ has a Newton polygon 
of type $[0,p-1,p]$.
To see that $v(\beta)=b/p$, observe that $X^c_\circ$ 
contains a point of $\cF(M)$, and $v(P'(\tau))=(p-1)v(\beta)$ is constant 
on $\cF(M)$ equal to $(p-1)b/p$. 
\hfill $\Box$

\begin{cor}\label{main-np-c}
Let $y\in Y^c(\Omega)$ and let $\alpha$ be a zero of the polynomial $g(y,X)$. 
Then $g(y,X+\alpha)$ has a Newton polygon of type $[0,p-1,p]$ 
and slope $b/p$ over $[0,p-1]$, 
where $e_{L/K}b=c(M/L)$ and $e_{L/K}$ is the ramification index of $L/K$.  
\end{cor}

Let $b$ be as in Theorem \ref{main-np-t}, 
$Z^{c,b}=Y^c\times D^{(b)}$,
$U=\psi^{-1}(Z^{c,b})$, and $u:U\rightarrow Z^{c,b}$ be the restriction 
of $\psi$.  
We put $Y^{c+}=\bigcup_{a>c} Y^a$, $D^{(b+)}=\bigcup_{a>b} D^{(a)}$,
$Z^{c+,b+}=Y^{c+}\times D^{(b+)}$ and $U^{+}=\psi^{-1}(Z^{c+,b+})$.

\begin{prop}\label{2finite}
We have 
\begin{itemize}
\item[{\rm a)}] $b\leq c$ and $X^c\subset U \subset X^b$;
\item[{\rm b)}] $\pi_0(U)\simeq \cF^c(L)\simeq \cF^c(M)$;
\item[{\rm c)}] the map $U^{+}\rightarrow Z^{c+,b+}$
is finite \'etale and geometrically totally decomposed,
and $\pi_0(U^+)\simeq\cF(M)$.
\end{itemize} 
\end{prop}
{\it Proof.}\ Let $y\in Y^c(\Omega)$,   
$g(y,.):D^1\rightarrow D^1$ be the map defined by the 
polynomial $g(y,X)$, and $\alpha\in \co$ be a zero of $g(y,X)$.
Corollary \ref{main-np-c} implies that 
$U_y=g(y,.)^{-1}(D^{(b)})=D(g(y,X))$ is geometrically connected 
and $g(y,.)^{-1}(D^{(b+)})$ is not geometrically connected.

a) By definition of $c$, the map $X^{c+}\rightarrow Y^{c+}\times D^{c+}$
is finite \'etale and geometrically totally decomposed. 
Therefore for any $y\in Y^{c+}(\Omega)$,
$g(y,.)^{-1}(D^{(c+)})$ is not geometrically connected.
We deduce that $b\leq c$ and $X^c\subset U\subset X^b$. 

b) The map $U\rightarrow Y^c\times D^{(b)}$ is finite flat, 
and the geometric fibers of the map $U\rightarrow Y^c$ 
are connected. Thus $\pi_0(U)\simeq \pi_0(Y^c)\simeq 
\cF^c(L)\simeq \cF^c(M)$. 

c) Let $Y^{c+}_\circ$ be a geometric connected component 
of $Y^{c+}$, $C^{c+}_\circ=\psi^{-1}(Y^{c+}_\circ\times 0)$, and $U^+_\circ=
\psi^{-1}(Y^{c+}_\circ\times D^{(b+)})$:
\[
\xymatrix{
C^{c+}_\circ \ar[r]\ar[rd]& U^{+}_\circ\ar[d]^{{\rm pr}^{+}}\\
& Y^{c+}_\circ}
\]
The map $C^{c+}_\circ\rightarrow Y^{c+}_\circ$ is finite 
\'etale and geometrically totally decomposed.
After a finite separable base change of $K$,  
let $s_1,\dots,s_p:Y^{c+}_\circ\rightarrow U^{+}_\circ$ be the sections
of ${\rm pr}^+$ obtained by splitting $C^{c+}_\circ$, and let $e=(0,1)\in D^n
\times D$.
We claim that for every $1\leq i\leq p$, the image of the open immersion
\[
\iota_i:
\begin{array}[t]{clcr}
Y^{c+}_\circ\times D^{(b/p +)}& \longrightarrow & D^n\times D\\
(y,z)& \longrightarrow & s_i(y)+ze
\end{array}
\]
is contained in $U^{+}_\circ$, and $U^+_\circ$ 
is the disjoint union of the images
of $\iota_i$. It is enough to check the claim above a point 
$y\in Y^{c+}_\circ(\Omega)$, then it follows from the fact that 
$g(y,X+\alpha)$
admits a Newton polygon of type $[0,p-1,p]$ and slope $b/p$. 
Moreover, the maps $\psi\circ \iota_i:Y^{c+}_\circ \times D^{(b/p+)}
\rightarrow Y_\circ^{c+}\times D^{(b+)}$ are isomorphisms. 
Thus $U^+\rightarrow 
Y^{c+}\times D^{(b+)}$ 
is finite \'etale and geometrically totally decomposed,
and the map $\pi_0(U^+)\rightarrow \pi_0(Y^{c+})$ is $p$ to $1$. 
The isomorphism $\pi_0(U^+)\simeq \cF(M)$ 
follows because $\pi_0(Y^{c+})\simeq \cF(L)$. 

\hfill $\Box$

\begin{prop}\label{g-fnp}
The rigid map $u:U\rightarrow Z^{c,b}$ is \'etale.
\end{prop}
{\it Proof.}\ Let $g'_X(Y,X)=\frac{d}{dX}g(Y,X)$, and let $(y,x)\in U(\Omega)$
be such that $g'_X(y,x)=0$. Since $U\subset X^b$, then $(y,x)$ defines a 
homomorphism of $\co_K$-algebras $\varphi: \co_M\rightarrow \co/\pi^b\co$,
and $y\in Y^c(\Omega)$ defines a homomorphism 
of $\co_K$-algebras $\psi: \co_L\rightarrow \co/\pi^c\co$
such that the following diagram is commutative 
\[
\xymatrix{\co_M\ar[r]^{\varphi}& \co/\pi^b\co\\
\co_L\ar[r]^{\psi}\ar[u]&\co/\pi^c\co\ar[u]}
\]
Let $\pi_M$ (resp. $\pi_L$) be a uniformizer of $\co_M$ (resp. $\co_L$).
Since $c>1$, then there is a unit $v\in \co^*$ such that
$\psi(\pi_L)=v\pi_L ({\rm mod} \ \pi^c\co)$
(see the proof of Proposition \ref{r-etale}). 
So $\varphi(\pi_L)=v\pi_L({\rm mod} \ \pi^b\co)$. 
We have $v(\pi_L)=1/e_{L/K}<b=c(M/L)/e_{L/K}$ because 
$M/L$ is not tamely ramified so $c(M/L)>1$. 
Let $e=e_{M/L}$ and $t\in \co_M^*$ be such that $\pi_M^e=t\pi_L$. 
It follows that $\varphi(\pi_M)^e=\varphi(t)v\pi_L({\rm mod} \ \pi^b\co)$. 
Thus $\varphi(\pi_M)=u\pi_M ({\rm mod} \ \pi^b\co)$ 
for a unit $u\in \co^*$. 

We have $\Omega^1_{\co_M/\co_L}\simeq \co_M/\pi_M^\alpha\co_M$, 
where $v(\pi_M^\alpha)=\alpha/e_{M/K}=(p-1)b/p$. We deduce that 
\[
\Omega^1_{\co_M/\co_L}\otimes_{\varphi}\co/\pi^b\co\simeq \co/\pi^{b(p-1)/p}
\co.
\]
On the other hand, we have 
\[
\Omega^1_{\co_M/\co_L}\otimes_{\varphi}\co/\pi^b\co\simeq
{\rm Coker}(\xymatrix{\co\ar[r]^{g'_X(y,x)}&\co}) \otimes \co/\pi^b
\co\simeq \co/\pi^b\co.
\]
We get a contradiction because $b>0$. \hfill $\Box$

\vspace{2mm}

Let $\cU$ and $\cZ^{c,b}=\fY^c\times \cD^{(b)}$ 
be the normalized integral models
of respectively $U$ and $Z^{c,b}$ over $\co$. 
By the finiteness theorem of Grauert 
and Remmert, the map $u$ extends to a finite 
map $u:\cU\rightarrow \cZ^{c,b}$. 

\begin{teo}\label{m-fnp}
The map $u:\cU\rightarrow \cZ^{c,b}$ is \'etale.
\end{teo}
{\it Proof.}\ We denote 
by $\spe:U\rightarrow \ocU$ and $\spe:Z^{c,b}\rightarrow 
\ocZ^{c,b}$ the specialization maps. We consider $\cF(M)\subset U(\Omega)$
and $\cF(L)\subset Z^{c,b}(\Omega)$. 
Then $\spe(\cF(L))\simeq \pi_0(Z^{c+,b+})$, $\spe(\cF(M))\simeq \pi_0(U^+)$,
and $u^{-1}(\spe(\cF(L)))=\spe(\cF(M))$.
We have $\pi_0(Z^{c+,b+})\simeq \cF(L)$, $\pi_0(U^+)\simeq \cF(M)$
and $U^+\rightarrow Z^{c+,b+}$ is finite \'etale and geometrically 
totally decomposed. We deduce that $u:\cU\rightarrow \cZ^{c,b}$ 
is \'etale above the points $\spe(\cF(L))$. Since the special fiber 
of $\cZ^{c,b}$ is smooth (by Theorem \ref{t-def})
and each of its geometric connected component contains 
a point of $\spe(\cF(L))$,
then $u:\cU\rightarrow \cZ^{c,b}$ is \'etale above the generic 
points of the special fiber of $\cZ^{c,b}$. 
Since $u:U\rightarrow  Z^{c,b}$ is \'etale (by Proposition \ref{g-fnp}), 
$\cU$ is normal, and $\cZ^{c,b}$ is regular (by Theorem \ref{t-def}), then 
by Zariski's purity theorem $u:\cU\rightarrow \cZ^{c,b}$ is \'etale. 
\hfill $\Box$

\vspace{2mm}

With the notation of the beginning of this section, we have 
\[
U=\Sp\left(\co_{Z^{c,b}}[x]/(P(x)-z)\right).
\]
Let $\tau$ be the residue class of $x$, which is a function in $\co_{\cU}$.

\begin{prop}
The polynomial $P(x+\tau)-z\in \co_{\cU}[x]$ has a Newton polygon of 
type $[0,p-1,p]$. Moreover, $P'(\tau)\co_\cU=\beta^{p-1}\co_{\cU}$
for an element $\beta\in \co$ with $v(\beta)=b/p$, where 
$e_{L/K}b=c(M/L)$ and $e_{L/K}$ is the ramification index of $L/K$.
\end{prop}
{\it Proof.}\ Similar to the proof of Theorem \ref{main-np-t} 
when we replace Proposition \ref{rec1} by Theorem \ref{m-fnp}.
\hfill $\Box$

\vspace{2mm}

The geometric monodromy induces an $\oOmega$-linear action 
of the inertia group $G^1$ on $\ocU$, $\ofY^c$ and $\ocD^{(b)}$. 

\begin{cor}\label{u-mon}
The group $G^c$ stabilizes the (geometric) 
connected components of $\ocU$, and the subgroup $G^{c+}$ acts trivially 
on $\ocU$.
\end{cor}
{\it Proof.}\ We know that $G^c$ stabilizes the geometric connected 
components of $X^c$. Thus it stabilizes those of $U$ because
$\pi_0(X^c)=\pi_0(U)$, and hence also those of $\ocU$.  
By Corollary \ref{d-c+} and Lemma \ref{d-mon}, $G^{c+}$ acts trivially 
on $\ofY^c\times \ocD^{(b)}$. The group $G^{c+}$ acts also 
trivially on $\spe(\cF(M))\simeq \pi_0(U^+)=\cF(M)$. 
Since the map $\ocU\rightarrow \ofY^c\times \ocD^{(b)}$ 
is \'etale and each connected component of $\ocU$ contains a point 
of $\spe(\cF(M))$, then $G^{c+}$ acts trivially on $\ocU$.
\hfill $\Box$

\vspace{2mm}

By Corollary \ref{main-np-c}, we have a well defined $K$-isomorphism
\[
\begin{array}{clcr}
U\times D^{(b/p)}&\longrightarrow & U\times_{Y^c}U\\
((y,x),z)&\longmapsto& ((y,x),(y,x+z))
\end{array}
\]
that makes $U$ as a $D^{(b/p)}$-torsor on $Y^c$. 
The normalized integral model of $U\times D^{(b/p)}$ over $\co$ 
is canonically isomorphic to $\cU\times \cD^{(b/p)}$. 
Theorem \ref{m-fnp} implies that the map $\cU\rightarrow \fY^c$ is smooth.
Then the 
normalized integral model of $U\times_{Y^c}U$ over $\co$ 
is canonically isomorphic 
to $\cU\times_{\fY^c}\cU$. Therefore, there is a natural
isomorphism 
\[
\cU\times \cD^{(b/p)} \longrightarrow  \cU\times_{\fY^c}\cU
\] 
that makes $\cU$ as a $\cD^{(b/p)}$-torsor on $\fY^c$. 
The induced isomorphism of special fibers 
\begin{equation}\label{torsor-mon}
\ocU\times \ocD^{(b/p)}\longrightarrow \ocU\times_{\ofY^c}\ocU
\end{equation}
is equivariant for the geometric monodromy action of $G^1$.

\section{The logarithmic ramification filtration}

\begin{prop} \label{relative}
Let $L$ be a finite separable extension 
of $K$, and $(Z,I,P)$ be a logarithmic system of generators 
of $\co_L$ over $\co_K$. 
Let $L'$ be a finite Galois extension 
of $L$. There exist 
\begin{itemize}
\item[{\rm (i)}] $(Z',I',P')$ a logarithmic system of generators 
of $\co_{L'}$ over $\co_K$ such that $Z'\supset Z$, so we can
identify $I$ with a subset of $I'$, and $P'\supset P$. We put $J=I'-I$;
\item[{\rm (ii)}] a rigid-analytic morphism 
$\varphi: D^{I'}\longrightarrow D^{I}\times_K D^{J}$, 
where the first factor is the canonical projection;
\end{itemize}
such that
\begin{itemize}
\item[{\rm a)}] $\varphi$
is finite and flat of degree $[L':L]$; 
\item[{\rm b)}] $\varphi^{-1}(\cF(L)\times \{0\})=\cF(L')$ 
and $\varphi$ is \'etale at every point in $\cF(L')$;
\item[{\rm c)}]  for any rational number $a>0$, 
there exist positive rational numbers $(a_j)_{j\in J}$ 
such that 
\[
Y^a_{Z',P'}= \varphi^{-1}(Y^a_{Z,P}\times\prod_{j\in J}D^{1,(a_j)}). 
\] 
\end{itemize}
\end{prop}
{\it Proof.}\ Let $L\subset L'\subset L''$ is a tower of finite Galois 
extensions of $L$. If we prove the Proposition for $L'/L$ and $L''/L'$,
then we deduce it for $L''/L$. Thus, we proceed by induction on $[L':L]$. 
We may assume that either $\co_{L'}/\co_L$ is \'etale, or 
$d=[L':L]$ is prime 
and the residue extension $\oL'/\oL$ 
is purely inseparable because the inertia group is solvable. 
In both cases, the extension $\co_{L'}/\co_L$ is monogenic 
generated by $z'\in \co_{L'}$.

First, we assume that the ramification index of $L'/L$ is $1$. 
Then $z'$ is a unit in $\co_{L'}$. We take $I'=I\amalg \{z'\}$ and $P'=P$. 
Let $F(X')\in \co_L[X']$ be the minimal polynomial of $z'$ over $\co_L$,
and $\tF\in \co_K[(X_i)_{i\in I}][X']$ be a monic polynomial 
of degree $[L':L]$ in $X'$ which lifts $F$. We easily see that 
\[
\begin{array}{clcr}
\varphi: &D^{I'}=D^I\times D^1&\longrightarrow & D^I\times D^1
\ \ \ \ \ \\
&\ \ \ \ \ \ \ \ \ (x,x')&\longmapsto& (x,\tF(x,x'))
\end{array}
\]
satisfies the required properties. Indeed, $Y^a_{Z',P'}
=\varphi^{-1}(Y^a_{Z,P}\times D^{1,(a)})$.

Second, we assume that $L'/L$ is totally ramified.
Then, we may assume that $z'$ is a uniformizer of $L'$. 
Let $\iota\in P$ be such that $z_\iota$ is a uniformizer of $L$. 
Since the minimal polynomial $F(X')$ of $z'$ over $\co_L$ is 
Eisenstein, then there exist $a_{d-1},\dots,a_0\in \co_L$ 
with $v(a_0)=0$, such that 
\[
F(X')=X'^d+z_\iota(a_{d-1}X'^{d-1}+\dots+a_0) \in \co_L[X']. 
\]
We take $I'=I\amalg \{z'\}$ and $P'=P\amalg\{z'\}$. 
For $0\leq i\leq d-1$, we fix $A_i\in \co_K[(X_i)_{i\in I}]$ 
a lifting of $a_i$. We put 
\[
\begin{array}{clcr}
\varphi: &D^{I'}=D^I\times D^1&\longrightarrow & D^I\times D^1
\ \ \ \ \ \ \ \ \ \ \ \ \ \ \ \ \\
&\ \ \ \ \ \ \ \ \ (x,x')&\longmapsto& (x,x'^d+x_\iota\sum_{i=0}^{d-1}
A_i(x)x'^i)
\end{array}
\]
Properties a) and b) are clearly satisfied. 
For c), we claim that 
\begin{equation}\label{eq134}
Y^a_{Z',P'}
=\varphi^{-1}(Y^a_{Z,P}\times D^{1,(a+1/e)})
\end{equation}
where $e$ be the ramification index of $L/K$. Let $e_i=v_L(z_i)$ 
(for $i\in I$), and put $X=(X_i)_{i\in I}$
and $(X,X')=(X_i)_{i\in I'}$. 
We fix liftings $g(X)$ of $z_\iota^e/\pi$ 
and $h_i(X)$ of $z_i/z_\iota^{e_i}$. 
The unit $z^d/z_\iota$ lifts to the polynomial
$h(X,X')=-\sum_{i=0}^{d-1}A_i(X)X'^i$. We choose a lifting 
$k(X,X')$ of the unit $z_\iota/z^d$. 
Then, $\tilde{g}(X,X')=g(X)h(X,X')^e$ is a lifting 
of $z^{de}/\pi$, and 
$\tilde{h}_i(X,X')=h_i(X) k(X,X')^{e_i}$ is a lifting of 
$z_i/z^{de_i}$ (for $i\in P$). 
We denote by $Y$ the right hand side of (\ref{eq134}). 
Obviously $Y\subset X^a_{Z'}$. We prove the other relations. 
Let $(x,x')\in Y(\Omega)$. 
We have $x'^d=x_\iota (h(x,x')+\pi^a *)$. Therefore, 
\begin{eqnarray*}
x'^{de}&=&x_\iota^e (h(x,x')+\pi^a *)^e=\pi(g(x)+\pi^a*)(h(x,x')+\pi^a *)^e\\
&=& \pi(g(x) h(x,x')^e+\pi^a*)=\pi \tilde{g}(x,x')+\pi^{a+1}*
\end{eqnarray*}
Since $h(X,X') k(X,X')$ lifts $1$, then $|h(x,x')k(x,x')-1|\leq \theta^a$. 
Observe that $v_{L'}(x')=1$. 
Thus, $x_\iota=x'^d(k(x,x')+\pi^a*)$. Then, for $i\in P$, we have
\begin{eqnarray*}
x_i&=& x_{\iota}^{e_i}(h_i(x)+\pi^a*)= x'^{de_i}(k(x,x')+\pi^a*)^{e_i}
(h_i(x)+\pi^a*)\\
&=& x'^{de_i}(k(x,x')^{e_i} h_i+\pi^a*)=x'^{de_i}\tilde{h}_i(x,x')
+\pi^{a+e_i/e}*
\end{eqnarray*}
It follows from (\ref{log-def2}) that 
$Y\subset Y^a_{Z',P'}$. The converse 
is obvious. \hfill $\Box$

\begin{cor}\label{absolute} 
Let $L$ be a finite Galois extension of $K$. 
There exist $(Z,I,P)$ a logarithmic system of generators of $\co_L$ 
over $\co_K$ and  
$f:D^I\rightarrow D^I$ a rigid-analytic morphism, such that the 
following properties hold~:
\begin{itemize}
\item[{\rm (i)}] $f$ is finite and flat of degree $[L:K]$;
\item[{\rm (ii)}] $f^{-1}(0)=\cF(L)$ and $f$ is \'etale above $0$;
\item[{\rm (iii)}] for any rational number $a>0$, there exist positive 
rational numbers $(a_i)_{i\in I}$ such that 
\[
Y^a_{Z,P}= f^{-1}(\prod_{i\in I}D^{1,(a_i)}).
\]
\end{itemize}
\end{cor}

\begin{prop} \label{log-section} 
Let $K\subset E \subset L$ be finite separable extensions. 
For a real number $a>0$, we put 
\[
\cF^{a-}_{\log}(L)=
\lim_{\stackrel{\longleftarrow}{0<b<a}}\cF^b_{\log}(L)\ \ \ {\rm and} \ \ \ 
\cF^{a+}_{\log}(L)=
\lim_{\stackrel{\longrightarrow}{b>a}}\cF^b_{\log}(L),
\]
where $b$ denotes a rational number. 
\begin{itemize}
\item[{\rm (i)}] For a rational number $a>0$, the map $\cF(L)\rightarrow 
\cF^a_{\log}(L)$ is surjective.  
\item[{\rm (ii)}] The map $\displaystyle \cF(L)\rightarrow 
\lim_{\stackrel{\longleftarrow}{a\in {\Bbb Q}_{>0}}}\cF^a_{\log}(L)$
is bijective. 
\item[{\rm (iii)}] $\cF^{a-}_{\log}(L)=\cF^{a}_{\log}(L)$  if $a$ is rational,
and $\cF^{a-}_{\log}(L)=\cF^{a+}_{\log}(L)$ if $a$ is not rational. 
\item[{\rm (iv)}] The following diagram is cocartesian
\[
\xymatrix{
\cF(L)\ar[r]\ar[d]& \cF^a_{\log}(L) \ar[d]\\
\cF(E)\ar[r]& \cF^a_{\log}(E)}
\]
\end{itemize}
\end{prop}
{\it Proof.}\ (i) Let $L'/K$ be a Galois closure of $L/K$. In 
the commutative diagram 
\[
\xymatrix{
\cF(L')\ar[r]^i\ar[d]_f&\cF^a_{\log}(L')\ar[d]^g\\
\cF(L)\ar[r]^j&\cF^a_{\log}(L)}
\]
$i$ is surjective by Corollary \ref{absolute}, 
and $g$ is surjective by Proposition \ref{relative}. 
Thus, $j$ is surjective.
(ii) follows from Proposition \ref{inj} for $r=+\infty$.
(iii) follows from Theorem \ref{cont} and (i).

We prove (iv). Let $L'$ be the Galois closure of $L/E$. 
It is enough to prove the Proposition after replacing $L/E$ by $L'/E$.
Hence, we may assume $L/E$ Galois. Then, the proof is similar to the proof of
Proposition \ref{section} (iii), if we use Proposition \ref{relative}
instead of Proposition \ref{finitude}. \hfill $\Box$

\begin{defi} Let $L$ be a finite separable extension of $K$, 
and $a>0$ be a rational number.
We say that the logarithmic ramification of $L/K$ is bounded 
by $a$ if the map $\cF(L)\rightarrow \cF^a_{\log}(L)$ is bijective. 
\end{defi}

\begin{prop}\label{log-fact-2}
The number 
\[
c_{\log}(L/K)=\inf\{ a\in {\Bbb Q}_{>0}\ ;\cF(L)\simeq \cF^a_{\log}(L)\}
\]
is rational, and the logarithmic ramification of $L/K$ is not bounded 
by $c_{\log}(L/K)$. 
\end{prop}
{\it Proof.}\ The same as Proposition 
\ref{fact-2} (i) and (ii). \hfill $\Box$

\subsection{Base change}

Let $K'$ be a finite separable extension of $K$ of ramification index $m$.
Let $L$ be a finite separable extension of $K$, and 
$L'$ be a finite separable extension of $K'$, such that $L\subset L'$.

\begin{lem}
Let $(Z,I,P)$ and $(Z',I',P')$ 
be logarithmic systems of generators of respectively $\co_L/\co_K$ 
and $\co_{L'}/\co_{K'}$. 
Let $a>0$ be a rational number, and let $Y^a_{Z,P}\subset D^I_K$ 
and $Y^a_{Z',P'}\subset D^{I'}_{K'}$ 
be the affinoid varieties associated respectively with
$\co_L/\co_K$ and $\co_{L'}/\co_{K'}$. 
Assume that $Z\subset Z'$, 
so we can identify $I$ with 
a subset of $I'$, and that $P\subset P'$. 
Then the canonical projection $D^{I'}_{K'}\rightarrow D^I_{K'}$
induces a rigid-analytic morphism
\begin{equation}\label{log-b-c-mor}
Y^{ma}_{Z',P'}\longrightarrow Y^{a}_{Z,P}\times_K K'. 
\end{equation}
\end{lem}
{\it Proof.}\ Let $x'=(x_i)_{i\in I'}\in Y^{ma}_{Z',P'}
(\Omega)$ and put 
$x=(x_i)_{i\in I}$. Obviously, we have $x\in X^a_Z(\Omega)$. 
We prove the other relations. Let $e$ be the ramification 
index of $L/K$, $e'$ be the ramification index of $L'/K'$,
and $r$ be the ramification index of $L'/L$. 
Let $e'_i=v_{L'}(z_i)$ (for $i\in P'$) and 
$e_i=v_L(z_i)$ (for $i\in P$). 
We have $re=me'$ and $re_i=e'_i$ for $i\in P$. 
We denote by $\pi'$ a uniformizer of $K'$, put 
$u=\pi'^m/\pi \in \co_{K'}^*$, and 
fix $\iota\in P'$  such that $z_{\iota}$ 
is a uniformizer of $L'$. 
We choose $g\in \co_{K'}[(X_i)_{i\in I'}]$ a lifting of 
$z_\iota^{e'}/\pi'$, 
$h_i \in \co_{K'}[(X_i)_{i\in I'}]$
a lifting of $z_i/z_\iota^{e'_i}$ (for $i\in P'$), and 
$k_i \in \co_{K'}[(X_i)_{i\in I'}]$
a lifting of $z_\iota^{e'_i}/z_i$ (for $i\in P'$).
Let $i\in P$.
We have $x_i=x_\iota^{re_i}(h_i(x')+\pi'^{ma}*)$. Hence, 
\begin{eqnarray*}
x_i^e&=& x_\iota^{ree_i}(h_i(x')+\pi^{a}*)^{e}
=x_\iota^{me'e_i}(h_i(x')+\pi^{a}*)^e\\
&=&\pi'^{me_i}(g(x')+\pi^{a}*)^{me_i}(h_i(x')+\pi^{a}*)^e
=\pi^{e_i} u^{e_i} g(x')^{m e_i} h_i(x')^e +\pi^{a+e_i}*
\end{eqnarray*}
Let $\bar{g}\in \co_K[(X_i)_{i\in I}]$ be a lifting of $z_i^e/\pi^{e_i}$.
Since 
\[
\frac{z_i^e}{\pi^{e_i}}=(\frac{z_i}{z_\iota^{re_i}})^e 
(\frac{z_\iota^{e'}}{\pi'})^{me_i} (\frac{\pi'^m}{\pi})^{e_i},
\]
then $u^{e_i} g^{me_i} h_i^e
\in \co_{K'}[(X_i)_{i\in I'}]$ lifts $z_i^e/\pi^{e_i}$.
Hence $|u^{e_i} g(x')^{me_i} h_i(x')^e-\bar{g}(x)|\leq \theta^a$. 
We deduce that $x_i^e=\pi^{e_i} \bar{g}(x) +\pi^{a+e_i}*$. 
Let $(i,j)\in P^2$. We have $x_j=x_\iota^{re_j}(h_j(x')+\pi^a*)$
and $x_\iota^{re_i}=x_i(k_i(x')+\pi^a*)$. Therefore, 
\[
x_j^{e_i}=x_\iota^{re_ie_j}(h_j(x')+\pi^a*)^{e_i}=x_i^{e_j}
(k_i(x')+\pi^a*)^{e_j}(h_j(x')+\pi^a*)^{e_i}=x_i^{e_j}
k_i(x')^{e_j}h_j(x')^{e_i}+\pi^{a+e_ie_j/e}*
\]
Let $\bar{h}_{i,j}\in \co_K[(X_j)_{j\in I}]$ be a lifting of $x_j^{e_i}
/x_i^{e_j}$.
Since 
\[
\frac{z_j^{e_i}}{z_i^{e_j}}=(\frac{z_j}{z_\iota^{re_j}})^{e_i}
(\frac{z_\iota^{re_i}}{z_i})^{e_j},
\]
then $h_j^{e_i}k_i^{e_j}$ lifts $x_j^{e_i}/x_i^{e_j}$. Hence,
$|h_j(x')^{e_i}k_i(x')^{e_j}-\bar{h}_{i,j}(x)|\leq \theta^a$. 
We deduce that $x_j^{e_i}=x_i^{e_j}\bar{h}_{i,j}(x)+\pi^{a+e_ie_j/e}*$. 
\hfill $\Box$

\begin{lem}\label{log-b-c}
We fix an embedding of $K'$ in $\Omega$, and assume that the 
natural map $\Hom_{K'}(L',\Omega)\rightarrow 
\Hom_{K}(L,\Omega)$ is injective. For a rational number 
$a>0$, if the logarithmic ramification of 
$\co_L/\co_K$ is bounded by $a$, then the logarithmic ramification of 
$\co_{L'}/\co_{K'}$ is bounded by $ma$.
\end{lem}
{\it Proof.}\ The same as Lemma \ref{b-c}. \hfill $\Box$

\begin{prop}\label{b-c-tame}
Assume that $K'/K$ is tamely ramified. 
Let $(L_j)_{j\in J}$ be finite separable extensions of $K'$ 
such that $L\otimes_K K'\simeq\prod_{j\in J} L_j$. The following 
are equivalent~:
\begin{itemize}
\item[{\rm (i)}] the logarithmic ramification of $\co_L/\co_K$ is 
bounded by $a$;
\item[{\rm (ii)}] the logarithmic ramification of 
$\co_{L_j}/\co_{K'}$ is bounded by $ma$ for one $j\in J$;
\item[{\rm (iii)}] the logarithmic ramification of 
$\co_{L_j}/\co_{K}$ is bounded by $a$ for one $j\in J$.
\end{itemize}
\end{prop}
We prove the proposition in several steps. 

\begin{lem}\label{red1}
Let $E$ be a finite separable extension of $L$ such 
that that $\co_{E}/\co_L$ is \'etale. 
Then Proposition $\ref{b-c-tame}$ holds for $L/K$ 
if and only if it holds for $E/K$. 
\end{lem}
{\it Proof.}\ 
Observe that the logarithmic ramification of $\co_{E}/\co_{K}$ 
is bounded by $a$ if and only if 
the logarithmic ramification of $\co_{L}/\co_{K}$ 
is bounded by $a$. One implication is obvious. 
Let $E_1$ be the maximal sub-extension of 
$E/K$ such that $\co_{E_1}/\co_{K}$ is \'etale.
Then $E$ is a composed extension of $L$ and $E_1$. 
Since the logarithmic ramification of $\co_{E_1}/\co_{K}$ is bounded 
by any rational number $b>0$, we get the other implication.
The Lemma follows by applying this equivalence twice 
(over $K$ and over $K'$). \hfill $\Box$

\begin{lem} \label{log-etale} 
Proposition $\ref{b-c-tame}$ holds if $\co_{K'}/\co_{K}$ is \'etale.
\end{lem}
{\it Proof}.\ (i)$\Rightarrow$(ii) follows 
from Lemma \ref{log-b-c}. (iii)$\Rightarrow$(i)
is obvious. It is enough to prove that 
(ii)$\Rightarrow$(iii) after replacing 
$K'$ by  a finite extension $K''$ such that 
$\co_{K''}/\co_K$ is \'etale. So, we may assume that 
$K'/K$ is Galois of group $H$. 
By Lemma \ref{red1}, we may assume 
that $L$ contains $K'$. 
Then $J=H$, and for $\sigma\in H$,  $L_\sigma=L$ equiped 
with the $K'$--algebra structure twisted  by $\sigma$.

We fix a monogenic presentation $\co_{K'}=\co_K[X]/P$, and 
denote by $z$ the image of $X$. Let $(Z=(z_i)_{i\in I},I,P)$ be a logarithmic
system of generators of $\co_L/\co_K$.
We may assume that there exists $\rho\in I-P$ such that $z_\rho=z$. 
Let $\co_L=\co_{K}[(X_i)_{i\in I}]/I_Z$ be the presentation associated 
with $Z$. We have 
\[
L_\sigma=\co_{K}[(X_i)_{i\in I}]/I_Z=
\co_{K'}[(X_i)_{i\in I}]/(I_Z, X_\rho-\sigma^{-1}(z))
\]
Let $Y^a$ and $Y_{\sigma}^a$ be the affinoid 
associated respectively with $(\co_{L}/\co_K,Z)$ and 
$(\co_{L_\sigma}/\co_{K'},Z)$. We put $Y=Y^a\times_KK'$. 
Then $Y^a_{\sigma}=Y(\theta^{-a}(X_\rho-\sigma^{-1}(z)))$. 
Let $x=(x_i)_{i\in I}\in Y(\Omega)$. 
Since $P(X_\rho)\in I_Z$, then 
\[
v(P(x_\rho))=\sum_{\sigma\in G}v(x_\rho-\sigma^{-1}(z))\geq a.
\]
We deduce that there exists a unique $\sigma\in G$ such that 
$v(x_\rho-\sigma^{-1}(z))\geq a$. For $\sigma'\not=\sigma$, 
$v(x_\rho-\sigma'^{-1}(z))=0$. It follows that the $Y^a_{\sigma}$,
for $\sigma\in H$, are disjoint and cover $Y$. 
To finish the proof, observe 
that $H$ acts naturally on $Y$, and 
permutes transitively the $Y^a_{\sigma}$ for $\sigma\in H$. 
\hfill $\Box$

\vspace{2mm}

\noindent {\it Proof of Proposition} \ref{b-c-tame}.\ 
There exist finite separable extensions $K_1/K$ and $K'_1/K'$
with $K_1\subset K'_1$~:
\[
\xymatrix{
K_1\ar[r]&K'_1\\
K\ar[u]\ar[r]& K'\ar[u]}
\]
such that~: a) $\co_{K_1}/\co_{K}$ and $\co_{K'_1}/\co_{K'}$ are \'etale;
b) $K'_1=K_1[X]/(X^m-\pi_1)$ where $\pi_1$ is a uniformizer of $K_1$;
and c) $\co_{K_1}$ contains a primitive $m$-th root of unity. 
By applying Lemma \ref{log-etale} to $K_1/K$ and to 
$K'_1/K'$, we are reduced to prove the Proposition under 
the assumption that $K'=K[X]/(X^m-\pi)$ and  
$\co_{K}$ contains a primitive $m$-th root of unity. 
It is enough to prove that (ii)$\Rightarrow$(iii). 
Let $\pi'\in \co_{K'}$ be
the image of $X$,
$\pi_L$ be a uniformizer of $L$ and $u=\pi/\pi_L^e$. By 
Lemma \ref{red1}, 
we may assume that there exists $v\in \co_L$ such that $v^m=u$. Therefore,
\[
\co_L\otimes_{\co_K}\co_{K'}=\co_{L}[X]/(X^m-\pi)=
\co_{L}[X]/(X^m-u\pi_L^e)\simeq \co_{L}[T]/(T^m-\pi_L^e)
\]
Let $d=(m,e)$, $m=dm'$ and $e=de'$. We have
\begin{eqnarray*}
\{\co_L\otimes_{\co_K}\co_{K'}\}'&\simeq& \prod_{\xi\in \mu_d(\co_K)}
\{\co_L[T]/(T^{m'}-\xi \pi_L^{e'})\}' \simeq 
\prod_{\xi\in \mu_d(\co_K)}
\{\co_L[T]/(T^{m'}- \pi_L^{e'})\}'\\
&\simeq& \prod_{\xi\in \mu_d(\co_K)}\co_L[T]/(T^{m'}- \pi_L)
\end{eqnarray*}
where $\{\}'$ denotes the normalization. 
For $\xi\in \mu_m(\co_K)$, we put $L_\xi=L[T]/(T^{m'}-\pi_L)$ equipped 
with the  morphism of $K$--algebras
\[
K'=K[X]/(X^m-\pi)\longrightarrow L[T]/(T^{m'}-\pi_L)
\]
given by $X\mapsto \xi v T^{e'}$. 
Obviously, 
the extensions $(L_\xi/K)_{\xi\in \mu_m(\co_K)}$ and $(L_j/K)_{j\in J}$
are all isomorphic. The set of isomorphism 
classes of the extensions $(L_\xi/K')_{\xi\in \mu_m(\co_K)}$ 
is equal to the set of isomorphism 
classes of the extensions $(L_j/K')_{j\in J}$. 

Let $(Z,I,P)$ be a logarithmic system of generators of $\co_L/\co_K$. 
We may assume that there exists $\iota\in P$ such that $\pi_L=z_\iota$. 
Let $\co_L=\co_K[(X_i)_{i\in I}]/I_Z$ be the presentation 
associated with $Z$. For $\xi\in \mu_{m}(\co_K)$, we have
\begin{eqnarray*}
\co_{L_\xi}&=&\co_K[(X_i)_{i\in I},T]/(I_Z,T^{m'}- X_\iota)\\
&=& \co_{K'}[(X_i)_{i\in I},T]/(I_Z,T^{m'}- X_\iota,\pi'- \xi vT^{e'})
\end{eqnarray*}
The image $t$ of $T$ in $\co_{L_\xi}$ is a uniformizer 
of $L_\xi$, and the ramification index of $L_\xi/K'$ is $e'$. 
We take $Z'=Z\amalg\{t\}$, $I'=I\amalg\{t\}$ and $P'=P\amalg\{t\}$.
Then $(Z',I',P')$ is a logarithmic system of generators for both
$\co_{L_\xi}/\co_K$ and $\co_{L_\xi}/\co_{K'}$. 
Let 
\[
\varphi_\xi: Y_\xi:=Y^{ma}_{Z',P'}\longrightarrow 
Y:=Y^a_{Z',P'}\times_K K'
\]
be the morphism (\ref{log-b-c-mor}) associated with the extensions
$L_\xi/K$ and $L_\xi/K'$. 
In the sequel, the variable $T$ is replaced by $X_t$.  Let 
$g\in \co_{K}[(X_i)_{i\in I}]$ be a lifting of $v\in \co_L^*$. 
By definition, there exist affinoid functions $f_1,\dots,f_r$ such that 
\begin{eqnarray*}
Y&=&D^{I'}(f_1,\dots,f_r,\theta^{-a-1}(\pi-g^m X_t^{me'}))\\
Y_\xi &=&D^{I'}(f_1,\dots,f_r,\theta^{-a-1/m}(\pi'-\xi g X_t^{e'}))
\end{eqnarray*}
Therefore, $Y_\xi=Y(\theta^{-a-1/m}(\pi'-\xi g X_t^{e'}))$. 
Let $x'=(x_i)_{i\in I'}\in \cY(\Omega)$ and put $x=(x_i)_{i\in I}$. 
For $\xi\not=\xi'$ in $\mu_m(\co_K)$,
we have $v((\xi-\xi') g(x) x_t^{e'})=1/m$. Since
\[
v(\pi-g(x)^mx_t^{me'})=\sum_{\xi\in \co_m(\co_K)}v(\pi'-\xi g(x) x_t^{e'})
\geq a+1,
\]
then, there exists a unique $\xi\in\mu_m(\co_K)$ such 
that $v(\pi'-\xi g(x) x_t^{e'})\geq a+1/m$. For $\xi'\not=\xi$, 
$v(\pi'-\xi' g(x) x_t^{e'})=1/m$. 
We deduce that the $Y_\xi$, for
$\xi\in \mu_m(\co_K)$, are disjoint and cover $Y$. 
To finish the proof, observe that ${\rm Gal}(K'/K)={\Bbb Z}/m{\Bbb Z}$
acts naturally on $Y$, and permutes transitively the $Y_\xi$
for $\xi\in \mu_m(\co_K)$. \hfill $\Box$

\subsection{Proofs of \ref{fil2}, \ref{log-properties} and \ref{log-rational}}
\label{ram-fil2}

Theorem \ref{fil2} Follows from Proposition \ref{log-section} (i), (ii)
and (iv). Theorem \ref{log-rational} follows from Proposition 
\ref{log-section} (iii). We prove now Proposition \ref{log-properties}. 
We omit the proof of the following easy Lemma. 

\begin{lem}\label{log-monogenic}
Let $L$ be a finite separable extension of $K$, 
$(Z,I,P)$ be a logarithmic system of generators  
of $\co_L/\co_K$ such that $\#P=1$, and $a>0$ be a rational number. 
Then, 
\begin{itemize}
\item[{\rm (i)}] $X^{a+1}_Z\subset Y^a_{Z,P}\subset X^a_Z$. 
\item[{\rm (ii)}] if moreover $Z=\{x\}$, and $x$ is a uniformizer of $\co_L$, 
then $Y^a_{Z,P}=X^{a+1}_Z$.
\end{itemize}
\end{lem}
Proposition \ref{log-properties}-3) is a consequence of Lemma 
\ref{log-b-c} and Proposition \ref{b-c-tame}.
By Lemma \ref{log-monogenic}-(ii), and unramified base change, 
Proposition \ref{log-properties}-4) reduces to Proposition \ref{properties}-3).
Proposition \ref{log-properties}-1) is a consequence 
of Lemma \ref{log-monogenic}-(i). 
We prove Proposition \ref{log-properties}-2). 
By Propositions \ref{properties}-1) and \ref{log-properties}-1),
we have $G_{\log}^{0+}\supset G^{1+}=P$, where $P$ is the wild inertia 
subgroup of $G$. For the converse, we consider a tamely ramified finite 
extension $L/K$. By applying Proposition \ref{b-c-tame} with $K'=L$, 
we deduce that $c_{\log}(L/K)=0$.
\hfill $\Box$

\appendix

\section{Eliminating fierce ramification}

Let $L$ be a finite separable extension 
of $K$, ${\frak d}_L\in \co_L$ be a different of $\co_L/\co_K$
(i.e. a generator of the different
defined in \cite{serre1} III), 
$\delta_L\in \co_K$ be a discriminant of $\co_L/\co_K$
(i.e. a generator of the discriminant defined in \cite{serre1} III), 
and $p^{s_L}$ be the inseparable degree of the residue extension $\oL/\oK$.

\begin{lem}\label{raynaud1}
Assume $s_L\not=0$, and let $x\in \oL$ be a radicial element over $\oK$. 
There exists a finite separable extension $K'$ of $K$ such that 
\begin{itemize}
\item[i)] $\pi=\pi_K$ is a uniformizing element of $\co_{K'}$;
\item[ii)] $\oK'=\oK[x]$; 
\item[iii)] for any composed extension $L'$ of $L$ and $K'$, 
the following inequality holds in the lexicographical order
\[
(v(\delta_{L'}),s_{L'})<(v(\delta_L),s_L),
\] 
where $\delta_{L'}\in \co_{K'}$ 
is a discriminant of $\co_{L'}/\co_{K'}$, and $p^{s_{L'}}$ is 
the inseparable degree of its residue extension. 
\end{itemize}
\end{lem}
{\it Proof.}\ Let $\overline{f}(X)=X^{p^e}-x^{p^e}\in \oK[X]$
be the minimal polynomial of $x$, and  $f\in \co_K[X]$ be a separable monic 
lifting of degree $p^e$. 
By \cite{serre1} I Corollary 1 of Proposition 15, $K'=K[X]/f$
is a finite separable field extension of $K$ with valuation ring 
$\co_{K'}=\co_K[X]/f$, and residue field $\oK[x]$.

Since $\oK'/\oK$ is purely inseparable, then 
$A=\co_L\otimes_{\co_K}\co_{K'}$ 
is a local ring. Let $L'$ be a composed extension of $L$ and $K'$. 
If $\co_{L'}\not=A$, then 
$v(\delta_{L'})<v(\delta_L)$ 
and we are done. If $\co_{L'}=A$, then the canonical  
surjective map 
\[
\oL\otimes_{\oK} \oK[x]=\oL[X]/(X^{p^e}-x^{p^e})=\oL[X]/(X-x)^{p^e}
\longrightarrow \oL'.
\]
shows that $\oL'\simeq \oL$. Therefore, $s_{L'}<s_L$. \hfill $\Box$

\begin{cor}\label{eliminating} 
There exist a finite 
separable extension $K'$ of $K$ and a composed extension 
$L'$ of $L$ and $K'$, such that 
\begin{itemize}
\item[i)] $\pi=\pi_K$ is a uniformizing element of $K'$;
\item[ii)] $L'/K'$ is unfiercely ramified.   
\end{itemize}
\end{cor}
{\it Proof}.\ 
An extension $L/K$ with $v(\delta_L)=0$ or $s_L=0$ is obviously 
unfierce. Since $v(\delta_L)$ and $s_L$ are non-negative integers,
the Corollary follows 
by applying Lemma \ref{raynaud1} finitely many times. 
\hfill $\Box$

\begin{prop}\label{car-tame}
For a finite separable extension $L/K$ 
of ramification index $e$, we have $v_L({\frak d}_L)\geq e-1$.
The equality holds if and only if $L/K$ is tamely ramified.  
\end{prop}
{\it Proof.}\  Let $J={\frak D}_{L}^{-1}$ be the largest 
fractional ideal of $L$ such that ${\rm Tr}_{L/K}(J)= \co_K$.
Then $J'=J{\frak m}_K$ is the largest fractional ideal such that 
${\rm Tr}_{L/K}(J')= {\frak m}_K$. 
Hence $I=J{\frak m}_K{\frak m}_L^{-1}$
is the smallest fractional ideal of $L$ satisfying 
${\rm Tr}_{L/K}(I)\supset \co_K$. 
Its inverse 
${\frak D}'_{L}={\frak d}_{L}{\frak m}_K^{-1}{\frak m}_L$
is an integral ideal. So we have the inequality $v_L({\frak d}_L)\geq e-1$. 
The equality is equivalent to that $\co_L$ is the smallest 
fractional ideal $I$ of $L$ such that ${\rm Tr}_{L/K}(I)=\co_K$. 
Since ${\rm Tr}_{L/K}({\frak m}_L)\subset {\frak m}_K$, the equality 
is further equivalent to that the map $\oL=\co_L/{\frak m}_L
\rightarrow \oK=\co_K/{\frak m}_K$ induced by ${\rm Tr}_{L/K}$ 
is non trivial. The induced map is equal to $e\times {\rm Tr}_{\oL/\oK}$
and is non-zero if and only if $e$ is prime to the residue characteristic
and $\oL$ is separable over $\oK$. \hfill $\Box$

\noindent CNRS UMR 7539, LAGA, Institut Galil\'ee, Universit\'e Paris-Nord, 
93430 Villetaneuse, France \\
E-mail: abbes@math.univ-paris13.fr

\vspace{3mm}

\noindent Department of Mathematics, University of Tokyo, Tokyo 113 Japan\\
E-mail: t-saito@ms.u-tokyo.ac.jp

\end{document}